\documentclass{article}
\usepackage[preprint]{colm2025_conference}

\usepackage{microtype}
\usepackage{graphicx}
\usepackage{subfigure}
\usepackage{booktabs} 

\newcommand{\bbR}{\mathbb{R}}
\newcommand{\bx}{\mathbf{x}}
\newcommand{\bz}{\mathbf{z}}

\usepackage{hyperref}

\definecolor{darkblue}{rgb}{0, 0, 0.5}
\hypersetup{colorlinks=true, citecolor=darkblue, linkcolor=darkblue, urlcolor=darkblue}

\usepackage{amsmath}
\usepackage{amssymb}
\usepackage{mathtools}
\usepackage{amsthm}
\usepackage{color}
\usepackage{multirow}
\usepackage{enumitem}
\usepackage[capitalize,noabbrev]{cleveref}

\theoremstyle{plain}

\theoremstyle{definition}

\theoremstyle{remark}

\title{Optimization Problem Solving Can Transition to Evolutionary Agentic Workflows}

\author{Wenhao Li\hspace{4pt}$^1$ \quad Bo Jin$^{1}$ \quad Mingyi Hong$^{3}$ \quad
Changhong Lu$^{2}$ \quad
Xiangfeng Wang$^2$ \\
$^1$Tongji University \quad $^2$East China Normal University \quad $^3$University of Minnesota
}

\begin{document}

\ifcolmsubmission
\linenumbers
\fi

\maketitle




\begin{abstract}
This position paper argues that optimization problem solving can transition from expert-dependent to evolutionary agentic workflows. 
Traditional optimization practices rely on human specialists for problem formulation, algorithm selection, and hyperparameter tuning, creating bottlenecks that impede industrial adoption of cutting-edge methods. 
We contend that an evolutionary agentic workflow, powered by FMs and evolutionary search, can autonomously navigate the optimization space, comprising problem, formulation, algorithm, and hyperparameter spaces. 
Through case studies in cloud resource scheduling and ADMM parameter adaptation, we demonstrate how this approach can bridge the gap between academic innovation and industrial implementation. 
Our position challenges the status quo of human-centric optimization workflows and advocates for a more scalable, adaptive approach to solving real-world optimization problems.
\end{abstract}

\section{Introduction}

Optimization methods are essential for solving real-world application problems as they enable efficient decision-making by maximizing benefits or minimizing costs within constraints~\citep{luenberger1984linear,boyd2004convex,chong2013introduction}. 
These techniques are widely used across industries, such as supply chain management~\citep{romero2000optimization}, manufacturing~\citep{dimopoulos2000recent}, and energy systems~\citep{dincer2017optimization}, to improve resource allocation, reduce waste, and enhance productivity. 
Additionally, optimization drives innovation in fields like machine learning~\citep{sun2019survey} and renewable energy~\citep{bernal2009simulation}, where it refines models and designs sustainable solutions. 
With advancements in computing power and algorithms, optimization has become a cornerstone for tackling complex challenges, improving efficiency, and fostering progress in diverse domains.

When addressing an application problem, which is typically articulated in natural language, the process of solving it generally involves the formulation of an optimization model (e.g., linear programming~\citep{dantzig2002linear}, mixed-integer programming~\citep{achterberg2013mixed}) that mathematically represents the problem's objectives, constraints, and decision variables. 
This model serves as a formal abstraction of the real-world scenario, enabling the application of rigorous analytical techniques.
Subsequently, an optimization algorithm or method (e.g., simplex method~\citep{nelder1965simplex}, branch and bound~\citep{boyd2007branch}, gradient descent~\citep{ruder2016overview}) must be designed or selected to solve the formulated model, or an existing optimization solver can be employed.

Both the optimization model and the solution method often involve hyperparameters, which are parameters that govern the structure of the model (e.g., regularization terms~\citep{schmidt2007fast}) or the behavior of the algorithm (e.g., step sizes~\citep{bazaraa1981choice}, convergence criteria~\citep{greenhalgh2000convergence}). 
These hyperparameters require careful tuning to ensure the model's accuracy and the algorithm's efficiency. 
Ultimately, through this structured approach, an optimal or near-optimal solution to the problem can be derived, providing actionable insights or decisions that align with the objectives. 
This systematic process underscores the interplay between problem formulation, algorithmic design, and parameter optimization in achieving practical solutions to real-world challenges.


Unfortunately, a noticeable gap exists between how optimization technologies are developed in academic research and how they are implemented in industrial practice~\citep{bixby2002solving,chen2023mind}. 
While academia often focuses on advancing theoretical frameworks and pushing algorithmic performance under idealized conditions, industry grapples with complex real-world constraints and legacy systems~\citep{reisman1994devolution,ormerod1997or,ormerod2002nature}. 
Many state-of-the-art methods proposed in research papers struggle with practical deployment due to computational scalability, implementation complexity, or incompatibility with existing infrastructure~\citep{bohanec1994trading,ibanez2018discretionary,sun2022predicting}. 
Conversely, industrial practitioners frequently rely on outdated, well-understood techniques because of concerns about reliability, maintainability, and the scarcity of specialized expertise needed to transfer academic insights into production~\citep{rossit2019visual,shin2021effects}.



This disconnect reflects deeper issues: 
Academic research typically strives for mathematical elegance and asymptotic performance guarantees, whereas industry emphasizes robust, maintainable solutions that meet rigid operational requirements~\citep{oliveira2016perspectives}. 
Bridging this gap has historically fallen to a small pool of academic experts, who must translate theoretical advances into workable solutions, creating bottlenecks that limit the scale and speed at which cutting-edge optimization can be adopted~\citep{van2010ordering}. 
Below, we highlight four specific but interconnected challenges that stem from this over-reliance on human expertise, and explain why a new, more agentic workflow is urgently needed:


\noindent\textbf{Exclusive Reliance on Specialists: Blind Spots in Real-World Optimization.}
Industrial problem formulations—and subsequent choices of algorithms and hyperparameter settings—often hinge on a small number of highly specialized individuals. 
While these experts can produce effective outcomes, their limited availability creates severe bottlenecks, and their subjective outlook can lead to overlooked use cases. 
Promising solutions remain theoretically trapped, seldom generalized, or rigorously tested outside of niche applications.


\noindent\textbf{Fragmented Innovation: When Academic Timelines Fail Industrial Urgency.}
The research community primarily measures success through theoretical breakthroughs (e.g., new convergence bounds or better worst-case complexity), which may not translate cleanly into industrial viability. 
This pace of academic evolution is ill-equipped to match rapid industrial changes in data, processes, and competitive requirements~\citep{simchi2014om}. 
As a result, companies struggle to adopt new methods that are seemingly out of sync with pressing operational realities, deepening the divide between theoretical potential and everyday practicality.


\noindent\textbf{Hyperparameter Volatility: The Achilles’ Heel of Practical Adoption.}
Misconfigured hyperparameters can undo many powerful optimization algorithms. 
Even slight deviations in learning rates, population sizes, or other settings can derail performance. 
Organizations without streamlined mechanisms for experimentation and automated tuning face tedious trial-and-error, often leading them to abandon novel algorithms in favor of stable, if suboptimal, legacy solutions. 
This hyperparameter fragility deters the very innovation that academic research strives to promote.


\noindent\textbf{The Engineering Gap: Why Cutting-Edge Methods Remain Stuck on Paper.}
Even when a promising academic algorithm exists, translating it into a reliable and resource-efficient industrial system can be prohibitively complex. 
Many breakthroughs hinge on specialized hardware, massive datasets, or pristine environments—conditions rarely found in production. 
Consequently, potential benefits remain unrealized: 
Companies continue relying on suboptimal workflows, while the latest academic advances gather dust in journals, seldom validated at scale.


\begin{figure*}[htb!]
    \centering
    \includegraphics[width=\linewidth]{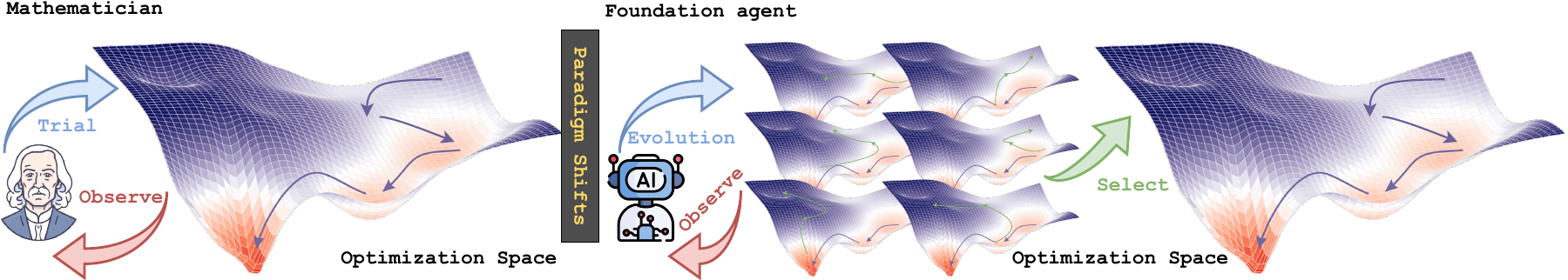}
    \caption{Comparative frameworks in optimization problem solving: human expertise and evolutionary agentic workflow.}
    \label{fig:paradigm-shifts}
\end{figure*}



These challenges illuminate how heavily modern optimization practice depends on a tight bottleneck of expert knowledge--one that struggles to keep up with dynamically shifting requirements and can stifle radical breakthroughs.
\textbf{This paper contends that an ``evolutionary agentic workflow,'' underpinned by foundation models (FMs), is uniquely suited to surmount these obstacles in ways that conventional optimization pipelines cannot.}
Rather than merely transferring a trendy technology from one domain to another, we argue that specific capabilities of FMs---including broad-domain knowledge integration, context-sensitive parsing of real-world constraints, and creative solution generation---align naturally with the evolving nature of industrial optimization problems.

FMs can parse incomplete or ambiguous descriptions of real-world challenges and translate them into coherent algorithmic or formulation strategies, drawing on extensive corpora of text, code, and domain data.
Concurrently, \textit{evolutionary search mechanisms}~\citep{yu2010introduction,zhou2011multiobjective,wu2024evolutionary} systematically explore vast spaces of potential configurations, iteratively testing and refining model structures and hyperparameters.
Such a workflow confronts core barriers in optimization \textit{because}:
\begin{enumerate}[leftmargin=*,label=(\alph*)]
    \item \textbf{Broad-Domain Generalization:} FMs cover diverse operational scenarios (e.g., supply chains, scheduling, manufacturing), bridging the domain fragmentation often faced when applying optimization theory to real-world constraints.
    \item \textbf{Iterative Refinement:} Evolutionary loops break away from one-shot human designs; they continuously discover new approaches, adapt parameters, and update problem formulations based on empirical feedback.
    \item \textbf{Creative Exploration vs. Rigorous Checking:} 
    Although large models may \textit{hallucinate}, sometimes proposing ideas that defy classical analysis, these ``illusions'' can seed genuinely novel heuristics if vetted by \textit{reasoning modules} (e.g., \texttt{OpenAI~o1}, \texttt{DeepSeek-R1}) or symbolic engines.
    \item \textbf{Dynamic Adaptation to Evolving Data:} Industrial realities shift rapidly, and FMs can swiftly reinterpret data, adapt constraints, and recommend altered strategies, mitigating reliance on static scripts.
\end{enumerate}

Beyond mere generation of text or code, we leverage specialized \textit{reasoning sub-components} to filter, refine, and logically validate proposals.
This reduces the risk that creative leaps devolve into mathematically unsound or operationally infeasible solutions.
In fact, \textbf{the agentic dimension} arises precisely because FMs alone, while rich in knowledge, do not inherently manage the end-to-end engineering pipeline:
\textit{reasoning modules} guide verification and symbolic checks, and \textit{evolutionary algorithms} orchestrate iterative solution refinement.
Together, these layers establish a feedback loop that continuously improves upon the search for new optimization ideas and the rigorous pruning of false leads.
Hence, we do not simply graft a AI paradigm onto optimization; we exploit FMs' specific, complementary strengths and evolutionary computing to tackle the persistent gaps identified above.


To arrive at this position, we introduce an \textit{optimization space} as the composition of (i) problem spaces, (ii) optimization formulation spaces, (iii) algorithm spaces, and (iv) hyperparameter spaces. 
This framework illuminates the essential questions of real-world optimization: 
\textit{``which problem formulation is appropriate?'', ``which algorithmic strategies are viable?'', and ``how should hyperparameters be tuned?''} 
Expert practitioners typically answer these questions leveraging domain knowledge, creating precisely the human bottlenecks we aim to alleviate.

We then propose anchoring agentic workflow~\citep{xi2023rise,cheng2024exploring} in heuristic search, particularly via evolutionary algorithms~\citep{yu2010introduction,zhou2011multiobjective,wu2024evolutionary} that mimic how human experts iteratively refine solutions in complex optimization landscapes. 
Evolutionary agentic workflow orchestrates this search by retaining and updating a dynamic memory of attempts (e.g., prior configurations and partial solutions), invoking specialized optimization tools (e.g., solvers) when needed, formulating plans for exploring or exploiting promising solution areas, and taking concrete actions to modify the current pipeline.
This framework creates a closed feedback loop: 
The agent's evolving understanding feeds back into subsequent search steps, enabling it to autonomously discover new algorithmic components, adapt hyperparameters, and innovate on problem formulations with minimal human intervention.

Moreover, we illustrate this evolutionary agentic workflow with two real-world case studies. 
In the first, the agent discovers novel heuristics for cloud resource scheduling—adjusting strategies to accommodate shifting workload patterns. 
In the second, it designs an adaptive step-size mechanism for ADMM (Alternating Direction Method of Multipliers), ensuring that the convergence properties remain stable. At the same time, performance is tuned to dynamic industrial conditions. 
Both studies prove that our hybrid of large language models and evolutionary search can realize meaningful gains in practical optimization tasks.
Finally, we address the current limitations of the agentic framework, particularly its lack of a built-in mechanism for theoretical verification and high inference cost.

\section{Optimization Space}

To systematically address the challenges outlined above and establish a framework for evolutionary agentic optimization workflows, we first need to formally characterize the space of optimization problems and their solution approaches. 
This formalization will serve as the foundation for understanding how foundation agents can navigate the complex landscape of real-world optimization tasks.
Formally, we can define an optimization space, which can be denoted as
\begin{equation}\label{eq:optimization-method-space}
    {\cal{O}} := {\cal{P}} \otimes {\cal{F}} \otimes {\cal{A}} \otimes {\cal{H}},
\end{equation}where ${\cal{P}}$ denotes the problem space, which contains the natural language or multi-modal description of application problems;
${\cal{F}}$ denotes the optimization formulation space, where each formulation can be represented using {\LaTeX} codes or other data structures (like an adjacency matrix for graph issues);
${\cal{A}}$ denotes the algorithm space, for instance, the heuristic methods for scheduling or planning, the gradient-type methods for continuous optimization, and the branch-and-bound techniques for mixed integer programming (MIP).
Each algorithm can also be represented by code without restricting the programming language.
{\cal{H}} denotes the hyperparameter space, which contains all the hyperparameters in the formulation and the algorithm, like the step-size (learning rate) and the regularization coefficient.

Each subspace can be a singleton.
For instance, we can design or search efficient algorithms for an optimization problem with a specific formulation like linear programming (LP), although the data can be freely changed.
In this case, the optimization formulation space contains only one specific formulation.
However, for the general case, all these sub-spaces can be an infinite-dimensional functional space.
Solving an optimization problem transforms into searching for the optimal solution in space ${\cal{O}}$ or the optimal combination of solutions in sub-spaces, respectively.

To illustrate how this formula captures the challenges of specialist dependency and hyperparameter sensitivity in real-world optimization, we examine two representative cases: cloud computing scheduling and distributed optimization. 
These examples, spanning discrete and continuous domains, will demonstrate why navigating the optimization space requires a more automated, agentic approach.

\noindent\textbf{Cloud Computing Scheduling.} 
In cloud computing, scheduling plays the key role, while a better scheduling policy can improve user experience and significantly save on purchasing computing resources~\citep{pietri2016mapping}.
Given a pool of computing resources that includes CPU, memory, GPU, etc., cloud computing providers aim to satisfy user computing requests through scheduling methods.
A scheduling problem $p_{schedule}\in {\cal{P}}_{schedule}$ is determined by the state of the resource pool ${\mathbf{s}}_{res}$ and the sequence of requests ${\mathbf{r}}_{req}$.
The scheduling problem $p_{schedule}$ is usually modeled as an online vector bin-packing problem, which can be represented by ${\mathbf{Bin}}$-${\mathbf{Packing}}$.
Heuristic methods (like best-fit or first-fit) are popularly designed to solve the online vector bin-packing model~\citep{bays1977comparison}.
Furthermore, combinatorial optimization~\citep{cote2021combinatorial} and machine learning~\citep{sheng2022learning} are also employed to solve this problem.
The algorithm space can be denoted as
$$
{\cal{A}}_{schedule}:=\big\{ {\mathbf{Heuristic}}, {\mathbf{CO}}, {\mathbf{ML}}, etc \big\}.
$$
Both in modeling and algorithm-design, there are some hyperparameters that need to be set, which can be denoted as ${\cal{H}}_{schedule}$.
Overall, following \eqref{eq:optimization-method-space}, solving the scheduling problem in cloud computing can be denoted as
\begin{equation}\label{eq:scheduling}
    {\cal{O}} := {\cal{P}}_{schedule} \otimes {\mathbf{Bin}}{\hbox{-}}{\mathbf{Packing}} \otimes {\cal{A}}_{schedule} \otimes {\cal{H}}_{schedule}.
\end{equation}

The challenge of navigating this optimization space lies in the unique characteristics of cloud workloads. 
For instance, when handling mixed CPU-GPU workloads, traditional bin-packing heuristics often fail to capture the complex resource dependencies and interference patterns between co-located tasks~\citep{wolke2015more}. 
While machine learning approaches can potentially learn these patterns~\citep{jiang2021learning,sheng2022learning}, they require careful formulation of the state space (e.g., whether to include historical resource utilization or job queue information) and reward signals (e.g., balancing immediate resource efficiency versus long-term cluster stability). 
The hyperparameter space further expands when considering practical constraints like power consumption limits and network topology, making manual tuning increasingly intractable.

\noindent\textbf{Distributed Optimization.}
ADMM~\citep{boyd2011distributed} has achieved significant success in distributed applications~\citep{yang2022survey}.
The classical ADMM algorithm is typically designed to solve the following structured convex optimization problems~\citep{boyd2011distributed}:
\begin{equation}\label{eq:convex}
\min_{\bx \in \bbR^n, \bz \in \bbR^m} f(\bx) + g(\bz),\quad \hbox{s.t.}\quad A \bx + B \bz = c,
\end{equation} where $f$ and $g$ are both proper, convex, and closed functions;
$A, B$ and $c$ establish the linear constraint.
A variety of optimization problems in machine learning~\cite{chang2020distributed}, statistics~\cite{boyd2011distributed}, and signal processing~\cite{hong2015unified} can be formulated in the form of \eqref{eq:convex}, and the application problem space can be denoted as ${\cal{P}}_{admm}$.
Besides those hyperparameters in function $f$ and $g$, there is a typical hyperparameter in ADMM method which is the penalty parameter in augmented Lagrangian function (usually denoted as $\beta$).
In practice, this parameter $\beta$ is sensitive and usually difficult to choose, as a result some self-adaptive schemes~\citep{he2016convergence} are proposed to dynamically tune $\beta$.
Overall, solving the structured optimization problem by ADMM can be summerized in the framework of \eqref{eq:optimization-method-space}, i.e.,
\begin{equation}\label{eq:admm-framework}
    {\cal{O}} := {\cal{P}}_{admm} \otimes \left\{ \eqref{eq:convex} \right\} \otimes {\mathbf{ADMM}} \otimes \left\{ \beta \right\}.
\end{equation}

The sensitivity to $\beta$ in ADMM reflects a deeper theoretical challenge: 
the penalty parameter affects both convergence speed and solution quality, but its optimal value depends on problem properties that are typically unknown a priori, such as the condition number of $A^TA$ and the Lipschitz constants of $f$ and $g$. 
While existing adaptive schemes like residual balancing~\citep{wohlberg2017admm} provide theoretical guarantees, they often perform poorly on ill-conditioned problems or when the objective functions have dramatically different scales. 
This suggests the need for more sophisticated adaptation strategies that can incorporate problem structure and runtime behavior.

\section{Evolutionary Agentic Workflow}

To effectively navigate the optimization space $\cal{O}$, we propose a novel framework that combines the reasoning capabilities of foundation agents with the systematic exploration power of evolutionary methods. 
Our framework consists of three key components: 
(1) foundation agents that leverage FMs for knowledge-driven optimization, 
(2) evolutionary methods that enable structured exploration of the solution space, and 
(3) an integrated evolutionary agentic workflow that orchestrates their interaction.

\begin{figure*}[htb!]
    \centering
    \includegraphics[width=\linewidth]{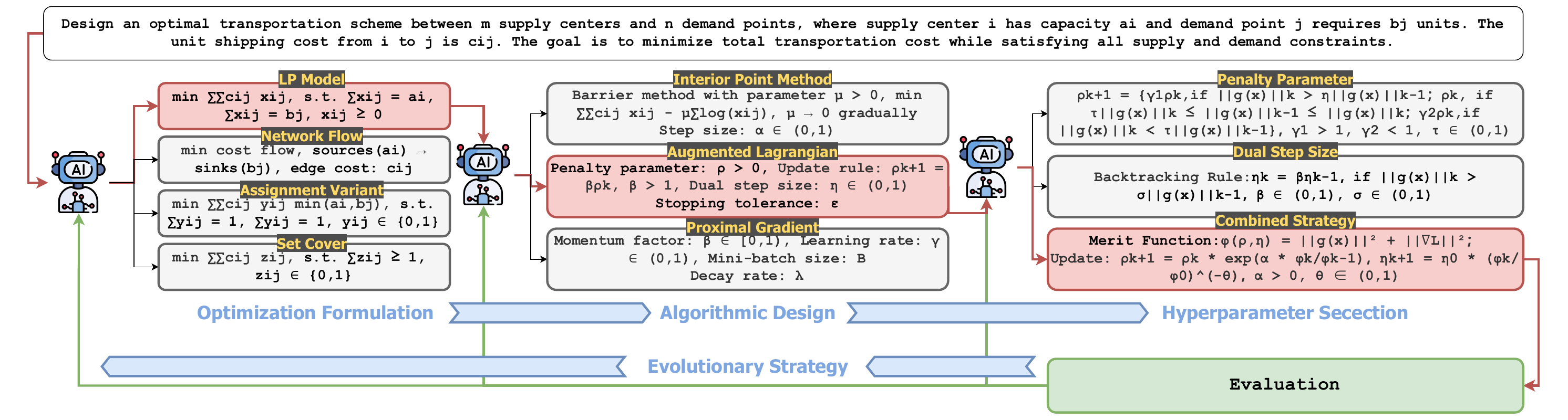}
    \caption{Automate optimization formulation, algorithmic design and hyperparameter selection for the transportation problem with evolutionary agentic workflow.}
    \label{fig:framework}
\end{figure*}

The foundation agents serve as intelligent guides in the optimization process, utilizing their broad knowledge base to make informed decisions about problem formulation, algorithm selection, and hyperparameter tuning. 
Meanwhile, evolutionary methods provide a systematic approach to explore and refine solutions, particularly valuable when dealing with large-scale or complex optimization problems where exhaustive search is impractical. 
The synergy between these components creates a powerful framework that balances exploitation of existing knowledge with exploration of novel solutions, as shown in Figure~\ref{fig:framework}.

\subsection{Foundation Agent Architecture}
\label{sec:foundation-agent-architecture}

We now describe the foundation agent~\citep{liu2024position}, which comprises $4$ interlinked modules:
\textit{Memory}, \textit{Reasoning}, \textit{World Modeling}, and \textit{Action}.
This design draws on the logic established in our introduction, where we emphasized the synergy between \emph{creative exploration} and \emph{rigorous checking}.
Specifically, the \emph{reasoning module} is tasked with maintaining mathematical and logical integrity, preventing the overall framework from being derailed by spurious model outputs or infeasible optimization proposals.

\subsubsection{Memory Module}
\noindent
The \textbf{memory module} enables foundation agents to store and retrieve knowledge and interaction histories that guide subsequent reasoning steps.
Similar to human memory~\citep{izquierdo1999separate}, it comprises:
(1) \textit{short-term memory} to capture the most recent context (such as the immediate conversation tokens or the latest hyperparameter guesses), and 
(2) \textit{long-term memory} to maintain cumulative encounters and insights gained across multiple optimization runs or experiments~\citep{wang2024voyager}.
Because the memory module can recall successful solution patterns, it helps the agent refine strategies based on previous attempts---a crucial mechanism for avoiding repeated mistakes~\citep{park2023generative,shinn2024reflexion}. 
For instance, if a particular linear solver consistently fails on specific constraints, the memory module can record this outcome, warning the agent against reusing that solver without modification.

\subsubsection{Reasoning Module}
\noindent
While FMs excel at wide-ranging text or code generation, they can inadvertently produce mathematically unsound or incomplete solutions, often called ``hallucinations.'' 
In order to reconcile \emph{creative exploration} with \emph{rigorous checking}, we introduce a dedicated \textbf{reasoning module} (implemented by large reasoning model~\citep{li2025system,xu2025towards,zhou2025large}, e.g., \texttt{OpenAI-o1}~\citep{openai2024o1} or \texttt{DeepSeek-R1}~\citep{guo2025deepseek}) that enhances chain-of-thoughts (CoT) depth and quality through the following technical route:
\begin{enumerate}[leftmargin=*,label=(\arabic*)]
    \item \textbf{Constructing High-Quality CoT Data for Optimization.}
    Unlike general-purpose tasks, real-world optimization often has complex constraints, domain-specific feasibility requirements, or shifting objectives.
    We compile CoTt datasets that reflect these subtleties: intermediate constraint checking, partial solution tracking, and sub-problem decomposition.  
    Such data directly encodes correct final solutions and the \emph{reasoning steps} necessary to confirm feasibility or assess trade-offs.

    \item \textbf{Domain-Focused Supervised Fine-Tuning (SFT).}
    We then fine-tune the base model on these specialized reasoning sequences,
    ensuring that it internalizes formal optimization logic (integer feasibility, objective bounding, dual formulation consistency, etc.).
    This fosters an informed CoT that can avoid naive heuristics and systematically explore realistic scenarios,
    aligning with the \emph{rigorous checking} agenda defined in our introduction.
    
    \item \textbf{Rule-Based Reward and Large-Scale Reinforcement Learning.}
    Finally, a rule-based reward mechanism evaluates each proposed CoTt:
    solutions that violate capacity constraints, fail integer feasibility, or produce inconsistent objective values are penalized,
    while improved or pragmatically validated solutions are rewarded.
    We then apply large-scale reinforcement learning to refine the model iteratively,
    boosting sequences that demonstrate methodical problem-solving and reducing reliance on ad hoc reasoning patterns. 
    Over successive cycles, the reasoning module converges toward producing disciplined, logically consistent solution paths
    that nonetheless retain an element of creativity in exploring new formulations.

\end{enumerate}

This reasoning framework is especially critical for large-scale optimization tasks that admit multiple local optima or require iterative constraint handling.
By integrating CoTt data rooted in \emph{actual optimization logic}, we nurture a system adept at diagnosing and correcting missteps early on.
In practice, the reasoning module can consult external symbolic or numeric engines to verify partial solutions,
fuse those validation signals into its CoTt, and then guide the agent’s next steps accordingly.

\subsubsection{World Modeling via LLMs}
\noindent
Prior research~\citep{hao2023reasoning,lin2024learning} has shown that large language models (LLMs) can serve as formidable \textbf{world models} by simulating the broad dynamics of a given system.
Here, an LLM-based world model can approximate the behavior of an optimization problem: 
it can hypothesize how constraints interact, estimate the cost of different solution branches, or predict how stochastic variables might evolve under certain assumptions~\citep{bai2022constitutional,ma2024eureka}.
Through CoTt~\citep{wei2022chain}, Tree-of-Thought~\citep{yao2024tree}, or related paradigms, an LLM can decompose a high-dimensional problem into manageable sub-problems~\citep{chu2023survey} and propose candidate pathways.
The ReAct framework~\citep{yao2022react} further equips agents with introspective checks to evaluate whether a reasoning path is consistent or sufficiently grounded. 
These generative and introspective aspects can be improved with environmental feedback~\citep{yao2022react,huang2022inner,wang2023describe,wang2024voyager,li2025multi}, 
FM-based critics~\citep{madaan2024self,chan2024chateval,shinn2024reflexion}, 
and human-in-the-loop guidance~\citep{huang2022inner}, complements the module’s latent knowledge.

\subsubsection{Action Module}
\noindent
Finally, the \textbf{action module} translates agent strategies into executable operations, bridging high-level plans and real implementation.
This may include the selection of specific problem formulations (linear vs.\ mixed-integer), 
algorithm choices (gradient descent vs.\ branch-and-bound),
and hyperparameter tuning (learning rates, tolerance levels, solver parameters).
Moreover, the action space can expand as external tools and specialized knowledge bases are integrated~\citep{lewis2020retrieval,schick2023toolformer,ge2023openagi,shen2024hugginggpt}.
By tapping into such resources, the action module can directly interface with modeling libraries or solver APIs, 
atlasing the current solution approach or pivoting to alternative methods if performance plateaus.

\medskip
\noindent
\textbf{Holistic Synergy in an Evolutionary Agentic Workflow.}\\
When these modules are combined within an evolutionary agentic workflow,
the agent alternates between exploratory variations of problem settings or parameter values and rigorous validation by the reasoning module.
Through ongoing memory updates, each promising direction accumulates historical evidence, 
informing the next generation of search steps guided by the reasoned CoTt.
Thus, the architecture balances the two imperatives highlighted in our introduction---%
\emph{creative exploration} and \emph{rigorous verification}---%
as a foundation for discovering powerful, real-world, viable optimization solutions.

\subsection{Evolutionary Framework}

While foundation agents excel at reasoning and knowledge utilization, they effectively exploring the vast optimization space $\cal{O}$ requires systematic exploration mechanisms.
Recent advances in combining LLMs with evolutionary approaches have shown promising results in various domains. 
Building upon these developments, we propose a generalized evolutionary framework that complements foundation agents in navigating the optimization space $\cal{O}$.
While traditional evolutionary methods rely on random mutations and predefined operators to generate new solutions, we propose a novel framework that leverages foundation agents as the solution generator within the evolutionary process. 
Our evolutionary framework incorporates three key mechanisms that guide the overall optimization process.

\noindent\textbf{Distributed Population Management.}
We maintain multiple solution populations (islands) that evolve independently, representing different exploration trajectories.
This distributed approach prevents premature convergence to local optima in $\cal{O}$ and provides natural checkpoints for solution diversity.
Within each island, foundation agents are called upon to generate new solutions by formulating specific problem models, designing concrete algorithms, or tuning particular hyperparameter combinations.

\noindent\textbf{Solution Diversity Preservation.} 
Solutions generated by foundation agents are clustered based on their characteristic signatures in the optimization space, which could include performance metrics, structural properties, or behavioral patterns.
This clustering mechanism helps maintain diversity across different dimensions of $\cal{O}$, preventing the search from collapsing to a single solution type.
The selection process balances between exploiting high-performing solutions and exploring diverse alternatives, informing the foundation agents about which solution characteristics to preserve or modify in subsequent generations.

\noindent\textbf{Knowledge-Guided Evolution.} 
When new solutions are needed, the evolutionary framework invokes foundation agents to generate them based on the current population state.
The agents leverage their memory and world models to understand successful patterns from existing solutions and create new variations.
This creates an efficient evolutionary cycle where high-performing solutions inform the agents' subsequent solution generation, while the evolutionary mechanism ensures systematic exploration.

\subsection{Human-Centered Evaluation}

Effective optimization requires algorithmic advancements and meaningful evaluation criteria that align with human priorities and domain-specific requirements. 
Traditional evolutionary methods often rely on predefined fitness functions that may not capture the nuanced objectives of real-world optimization problems. 
Drawing inspiration from human-centered AutoML approaches~\citep{lindauer2024position}, the proposed framework incorporates evaluation mechanisms that explicitly account for human preferences and domain expertise. 
Rather than optimizing solely for predictive performance, the evaluation considers multiple objectives, including interpretability, computational efficiency, and domain-specific constraints. 
The system enables domain experts to express preferences between candidate solutions through pairwise comparisons, effectively learning implicit utility functions that would be difficult to formalize. 

This preference-based approach~\citep{giovanelli2024interactive} allows the evolutionary process to navigate toward solutions that satisfy formal metrics and tacit domain knowledge. 
Additionally, the framework supports interactive refinement of the evaluation criteria as the optimization progresses, recognizing that human stakeholders often clarify their requirements through an iterative process of exploring the solution space. 
By establishing this bidirectional feedback loop between automated search and human evaluation, the framework can address the ``exclusive reliance on specialists'' challenge identified earlier, democratizing access to optimization while preserving the crucial role of human expertise in guiding the search toward meaningful solutions.

\noindent\textbf{Remark: Complementarity with AutoML.}
It is important to clarify how the evolutionary agentic workflow relates to existing AutoML frameworks. 
While AutoML systems excel at efficiently searching predetermined configuration spaces with well-defined boundaries, evolutionary agentic workflow addresses fundamentally different challenges in optimization. 
Unlike AutoML platforms that primarily focus on selecting and tuning algorithms within fixed modeling paradigms, evolutionary agentic workflows can navigate unclear formulation spaces where the very structure of the problem representation may be reconsidered. 
For instance, the framework can autonomously reformulate cost functions with different constraints or decision variables—transformations that traditionally require domain experts' manual intervention. 
Additionally, while AutoML excels at algorithm selection from established options (e.g., choosing between random forests and neural networks), it is not designed to auto-design entirely new heuristics or derive specialized algorithmic expansions like the ADMM variants described in our case studies.

We view evolutionary agentic workflow as complementary to, rather than competing with, traditional AutoML. 
Indeed, evolutionary agentic workflows can incorporate standard HPO (hyperparameter optimization) modules (e.g., Optuna~\citep{akiba2019optuna}) as components within their iterative search process, unifying the strengths of both paradigms. 
This integration allows evolutionary agentic workflows to leverage the efficiency of AutoML for well-defined subproblems while maintaining the creative exploration capabilities needed for less bounded optimization spaces. 
Furthermore, the framework aligns with recent calls for more human-centered AutoML approaches~\citep{lindauer2024position} by supporting flexible human feedback mechanisms through our evaluation module. 
This allows domain experts to provide targeted guidance at critical junctures while the system autonomously handles repetitive low-level tasks. 
This balanced division of labor creates a genuinely collaborative environment where machine creativity and human expertise contribute to discovering novel strategies that remain undiscovered through purely automated or manual approaches.

\subsection{Operational Workflow}

Integrating evolutionary methods and foundation agents creates a robust optimization framework where evolutionary search guides the overall exploration while foundation agents serve as sophisticated solution generators.
The evolutionary process begins with an initial population of solutions, each generated by foundation agents based on their knowledge of the problem domain.
As the evolution proceeds, the framework identifies promising solution characteristics and directs the agents to generate new solutions that build upon these successful patterns.

Foundation agents are called upon during each iteration to perform specific tasks within the broader evolutionary search.
They analyze existing solutions to understand successful patterns, generate new solution variations through problem reformulation or algorithm modification, and provide quality assessments of newly generated solutions.
The framework then manages these solutions across multiple islands, maintains population diversity, and determines which solutions should inform the next generation.

This structured workflow leverages the complementary strengths of both components: 
The evolutionary framework provides systematic exploration and selection pressure, while foundation agents contribute domain knowledge and sophisticated solution generation capabilities.
The resulting system demonstrates remarkable effectiveness in navigating complex optimization landscapes, as the evolutionary process efficiently explores the solution space. 
Meantime, foundation agents ensure each generated solution incorporates domain expertise and learning from previous successes.

\subsection{Generalization and Lifelong Learning}

A critical consideration for agentic optimization workflows is their potential ability to handle problems that extend beyond the foundation model's training distribution. 
While FMs demonstrate impressive zero-shot capabilities across many domains, real-world optimization often involves novel constraints, emerging domains, or unprecedented problem structures. 
Here, we discuss several promising mechanisms that could enhance generalization and support continual adaptation in evolutionary agentic workflows.

Retrieval-Augmented Generation (RAG) represents a compelling approach to extend the knowledge boundaries of FMs. 
By connecting the foundation agent to external knowledge repositories containing domain-specific papers, technical reports, and implementation examples, RAG could effectively enlarge the agent's knowledge scope beyond its pre-trained parameters. 
When encountering a problem with unfamiliar characteristics, such a system could extract key features from the problem statement, query the external knowledge base for relevant optimization techniques, and integrate retrieved information with its parametric knowledge to propose appropriate formulations or algorithms~\citep{lewis2020retrieval}. 
This hybrid approach would balance the foundation model's general reasoning capabilities with precise, specialized knowledge required for domain-specific optimization challenges.

For domains requiring more profound expertise, targeted supervised or reinforcement fine-tuning could significantly enhance an agent's ability to navigate specialized optimization landscapes~\citep{chen2023fireact,liao2025marft}. 
Even a modest number of examples (typically $10$-$50$ samples) might substantially improve performance on domain-specific tasks, particularly when these examples demonstrate common constraint patterns, typical failure modes, and their remediation, and effective problem decomposition strategies~\citep{bai2022constitutional,rft}. 
This approach would be particularly valuable when adapting the workflow to new industries or technical domains where problem structures differ substantially from standard benchmarks.

Continual (lifelong) learning mechanisms represent perhaps the most promising direction for ensuring long-term viability of agentic optimization workflows~\citep{zheng2025towards,zheng2025lifelong}. 
Such mechanisms would allow the optimization agent to evolve with experience as it encounters and solves new problems~\citep{silver2025experience}. 
This capability would maintain relevance in dynamic industrial environments where problem specifications continuously evolve. 
Potential continual learning processes could include experience replay, where the agent maintains a buffer of previously solved problems; 
selective parameter updates for critical performance improvements; 
and meta-learning approaches, where the agent could learn to adapt quickly to new problem classes by identifying common patterns across different optimization tasks.

This multi-faceted approach to generalization and lifelong learning could ensure that evolutionary agentic workflows remain effective even as optimization problems evolve beyond their original specifications. 
Rather than providing a static solution, such a framework would continue to adapt and improve through ongoing interactions with the optimization environment, mirroring how human experts develop and refine their problem-solving strategies over time.

\subsection{Component-wise Optimization}

While our framework proposes an integrated evolutionary approach where foundation agents serve as solution generators, examining how foundation agents have been applied to individual optimization components is instructive. 
Though not yet incorporating evolutionary mechanisms, current research has shown promising results in using foundation agents for problem formulation, algorithm design, and hyperparameter tuning separately~\citep{liu2024systematic}. 
These works provide valuable insights into the foundation agents' capability as solution generators and highlight the benefits of incorporating evolutionary methods. 

\noindent\textbf{Optimization Formulation.}
Recent works have demonstrated the potential of foundation agents in converting natural language descriptions into mathematical optimization models~\citep{zhang2019gap,meadows2022survey,wu2022autoformalization,wasserkrug2024large}. These efforts make optimization more accessible to non-experts by automating the problem formulation process.

Several frameworks have been proposed to leverage LLMs for optimization modeling. The NL4Opt competition~\citep{ramamonjison2023nl4opt}, inspired by OptGen~\citep{ramamonjison2022augmenting}, catalyzed significant progress in this direction, with various approaches achieving promising results~\citep{ning2023novel,gangwar2023highlighting}. Subsequent works like OptiMUS~\citep{ahmaditeshnizi2023optimus}, ORLM~\citep{tang2024orlm}, LM4OPT~\citep{ahmed2024lm4opt}, and MAMO~\citep{huang2024mamo} have further advanced the field by developing specialized models and techniques for optimization problem formulation.

Different approaches have been explored to enhance the formulation process. 
Some works focus on decomposing the modeling into subtasks~\citep{tsouros2023holy,khot2023decomposed}, while others emphasize interactive refinement with users~\citep{mostajabdaveh2024optimization,almonacid2023towards}. 
Various evaluation methods have been proposed, from comparing with ground truths~\citep{amarasinghe2023ai} to expert validation~\citep{li2023synthesizing,chen2024diagnosing}. 
Recent studies~\citep{fan2024artificial} have shown that while LLMs excel at textbook-level problems, they may require iterative refinement for complex real-world scenarios.

While these works demonstrate the potential of foundation agents in problem formulation, they primarily focus on generating single formulations rather than systematically exploring the vast space of possible problem models. 
Evolutionary methods could enhance these approaches by enabling structured exploration and progressive refinement.

\noindent\textbf{Algorithmic Design.}
Researchers have explored two main approaches to leverage foundation agents~\citep{huang2024large}. 
The first approach uses LLMs directly as optimization algorithms through prompting. 
For instance, OPRO~\citep{yang2024large} guides LLMs to generate solutions through in-context learning, incorporating previous solutions to enable iterative improvement. 
However, such direct approaches often face limitations when dealing with complex problems involving large search spaces~\citep{zhang2024understanding}, leading researchers to explore more sophisticated methods like utilizing LLMs to generate solver scripts~\citep{ahmaditeshnizioptimus,huang2024mamo}.

The second approach employs foundation agents to design optimization algorithms through agentic workflows. 
Some pioneering works have already demonstrated the potential of combining agentic workflows with evolutionary methods. 
FunSearch~\citep{Funsearch} and Evolution of Heuristics~\citep{EoH} integrate LLMs with evolutionary computation to automatically design heuristic algorithms, showing superior performance over manual designs. 
Building upon this direction, works like ReEvo~\citep{ye2024reevo} and AEL~\citep{liu2023algorithm,liu2024example} further enhance the evolutionary framework with reflection mechanisms and automated algorithm design capabilities.

Various specializations within this approach have emerged, from animal-inspired metaheuristics~\citep{zhong2024leveraging} to language model-based crossover operators~\citep{meyerson2024language}. 
Some researchers focus on population-based algorithms~\citep{pluhacek2023leveraging}, while others explore multi-objective optimization of algorithm design~\citep{yao2024multi}. 
These works demonstrate that evolutionary agentic workflows can effectively navigate the complex space of algorithm designs.

While these studies validate our position on the synergy between evolutionary methods and foundation agents, they have focused on applying this combination to algorithmic design in isolation, rather than integrating it into a comprehensive optimization framework.

\noindent\textbf{Hyperparameter Selection.}
While hyperparameter optimization is crucial for optimization algorithms, research on using foundation agents for this purpose has primarily focused on machine learning applications rather than optimization problems. 
Recent works have demonstrated promising approaches using agentic workflows for automated hyperparameter tuning~\citep{liu2024large,zhang2023using,mahammadli2024sequential}.

Despite these advances in ML, the application of foundation agents for hyperparameter selection in optimization algorithms remains relatively unexplored. 
Moreover, unlike the evolutionary approaches in algorithmic design, current hyperparameter optimization methods using foundation agents have not yet systematically incorporated evolutionary mechanisms to explore the hyperparameter space.




\section{Case Studies}

We present two case studies to validate the potential of evolutionary agentic workflows in optimization.
For the computational overhead, each workflow iteration takes about $2$ minutes and totaling $200$-$300$ rounds in agentic VM scheduling, and about $10$ seconds with $15$-$20$ total iterations in agentic ADMM step-size tuning.
These computations rely on calling external LLM APIs. 
Thus, local resource usage remains minimal at present. 
However, we anticipate that techniques like quantization, distillation, and pruning—as well as hardware-accelerated GPU kernels~\citep{lange2025ai,chen2025automating}—will reduce runtime if models are fully deployed on-site.

\subsection{Agentic Virtual Machine Scheduling}

Our first case study examines algorithm design for virtual machine scheduling in cloud computing, where we fix the problem formulation and hyperparameter settings to isolate the algorithm design space. 
This represents a special case of our framework where the optimization spaces for problem formulation and hyperparameter tuning are reduced to singletons, allowing us to focus solely on evolving scheduling algorithms through agentic workflows.

\begin{figure}[htb!]
    \centering
    \includegraphics[width=.7\linewidth]{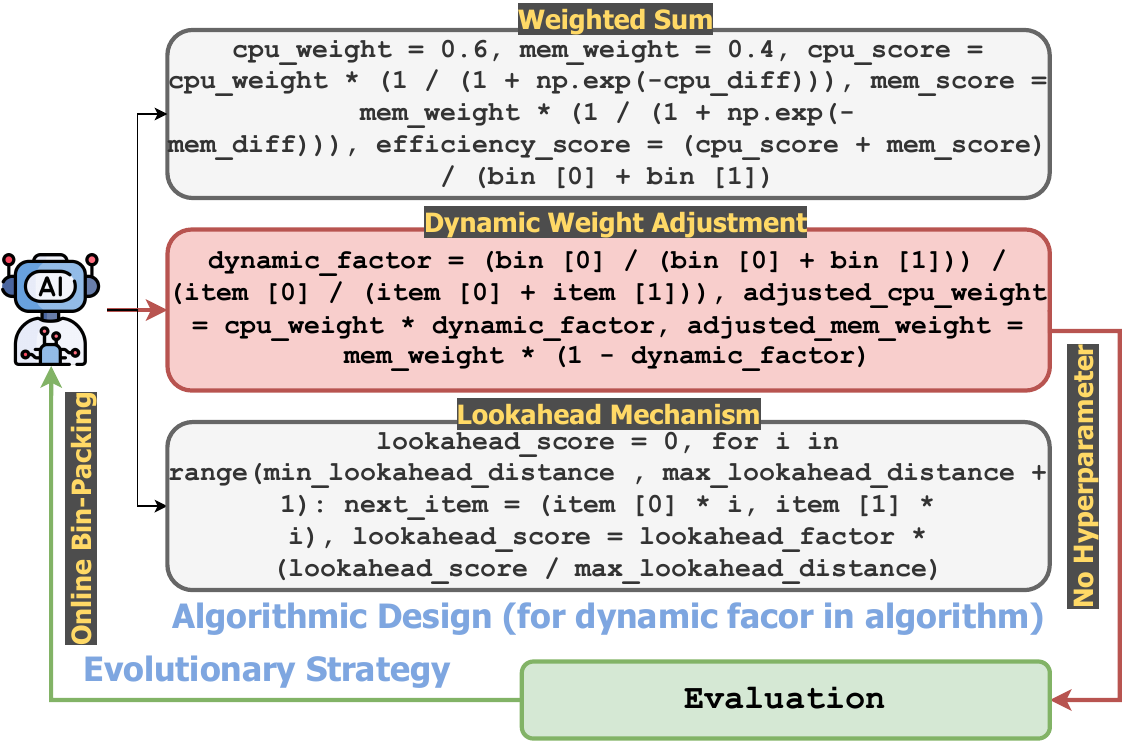}
    \caption{Evolutionary agentic workflow for VM scheduling.}
    \label{fig:vms}
\end{figure}

The virtual machine scheduling aims to design a scheduler to arrange the requests received by the proper virtual machines in the cloud computing system.
As discussed in the introduction, this application problem can be modeled as an online vector bin-packing problem, which is significantly challenging to solve.
The system needs to arrange the industry requests as soon as possible.
Heuristic methods have become the key to the online scheduling problem.
However, given a particular scheduling algorithm, it cannot always be efficient even for the request sequence from a relatively robust scenario.
This further motivates us to evolve the scheduling algorithm scheme, especially following the agentic workflow, as shown in Figure~\ref{fig:vms}.



We have conducted some preliminary experiments.
The datasets used are the traces collected by the Huawei Cloud~\footnote{\url{https://github.com/huaweicloud/VM-placement-dataset}.}
Each VM request in the dataset needs two types of resources: CPU cores and memory. 
For the baseline selection, we compare with a traditional heuristic method {\em{BestFit}}, which picks the server with the highest current allocation rate, aiming for an optimal use of resources.
We directly implement the procedure on the scenario with $50$ virtual machines, with the Deepseek coder LLM model~\footnote{\url{https://deepseekcoder.github.io}.}
One A100 (80G) GPU is employed for training, with $300$ training epochs and about $7$ to $8$ hours.
Furthermore, more scenarios with virtual machine size $30$, $100$, $150$ and $200$ are tested to show not only efficiency but also generalization ability.
The scheduling length is used to evaluate the performance, which is defined as the total number of VMs processed by the cluster over a given capacity.
The longer scheduling length indicates a more effective algorithm, reflecting the ability to manage and process more tasks over time.
In Table \ref{table:VM-scheduling}, we can find that the new algorithm performs significantly better than BestFit on the training scenario.
Furthermore, the new algorithm can beat BestFit if we directly implement it in other scenarios without training.

\begin{table}[!ht]
    \centering
    \resizebox{.6\linewidth}{!}{%
    \begin{tabular}{c|c|c|c}
    \hline
        Scenarios & VM Size & BestFit & \textbf{\textcolor{purple}{NewAlgorithm}} \\ \hline
        Train & 50 & 1386.8 & \textbf{\textcolor{purple}{1526.1}} \\ \hline
        Test & 30 & 566.9 & \textbf{\textcolor{purple}{572.6}} \\ \hline
        Test & 100 & 3811.5 & \textbf{\textcolor{purple}{3874.9}} \\ \hline
        Test & 150 & 6256.2 & \textbf{\textcolor{purple}{6318.2}} \\ \hline
        Test & 200 & 8366.5 & \textbf{\textcolor{purple}{8447.9}} \\ \hline
    \end{tabular}%
    }
    \label{table:VM-scheduling}
    \caption{Comparision between the NewAlgorithm and BestFit.}
\end{table}

\begin{figure}[htb!]
    \centering
    \includegraphics[width=0.5\linewidth]{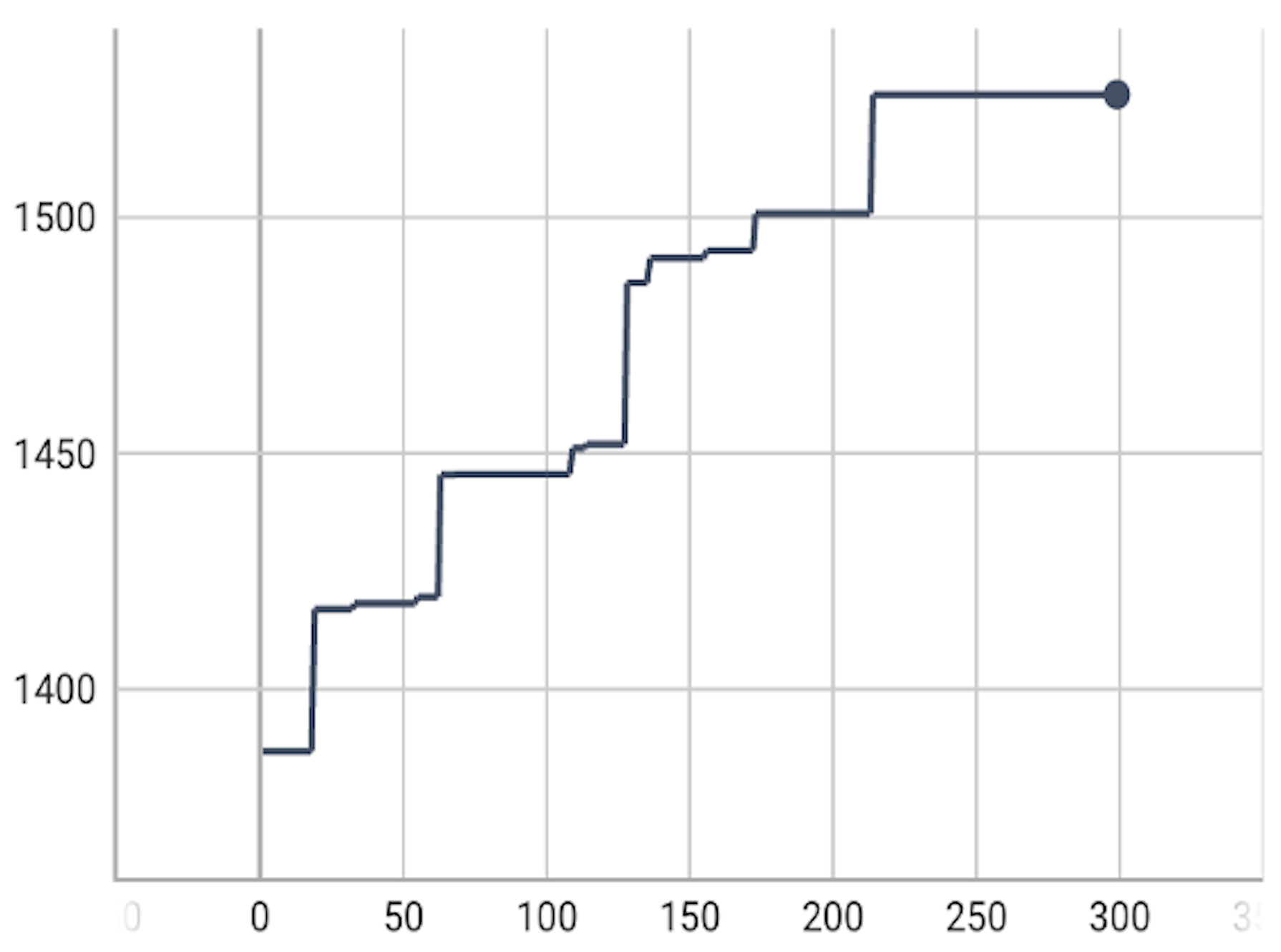}
    \caption{The perforcement of the discovered algorithm during traing on scenario with VM size $50$.}
    \label{fig:vm-scheduling}
    \end{figure}

We also demonstrate the evolution curve of algorithmic capabilities during the training process.
It can be seen in Figure \ref{fig:vm-scheduling} that the evolutionary generation algorithm is fully automated and its capabilities are progressively enhanced, significantly improving the efficiency of algorithm design.

\subsection{Agentic ADMM Step-Size Tuning}

Our second case study focuses on hyperparameter tuning for the ADMM, where we consider a fixed problem formulation and algorithm design. 
This represents another special case of our framework, where the optimization spaces for problem formulation and algorithm design are singletons, allowing us to concentrate on evolving hyperparameter selection strategies through agentic workflows. 
This fascinating case demonstrates how our framework can enhance even well-established optimization algorithms through adaptive parameter tuning.


The ADMM has shown significant success in various optimization problems in ML, statistics, and signal processing.
However, its performance is sensitive to the selection of the hyperparameter introduced in the augmented Lagrangian function (typically denoted as $\beta$), which is usually fixed in conventional implementations.
As discussed in the introduction section, employing ADMM to solve a structured convex optimization problem can be formulated in a unified framework as \eqref{eq:admm-framework}.

\begin{figure}[htb!]
    \centering
    \includegraphics[width=.7\linewidth]{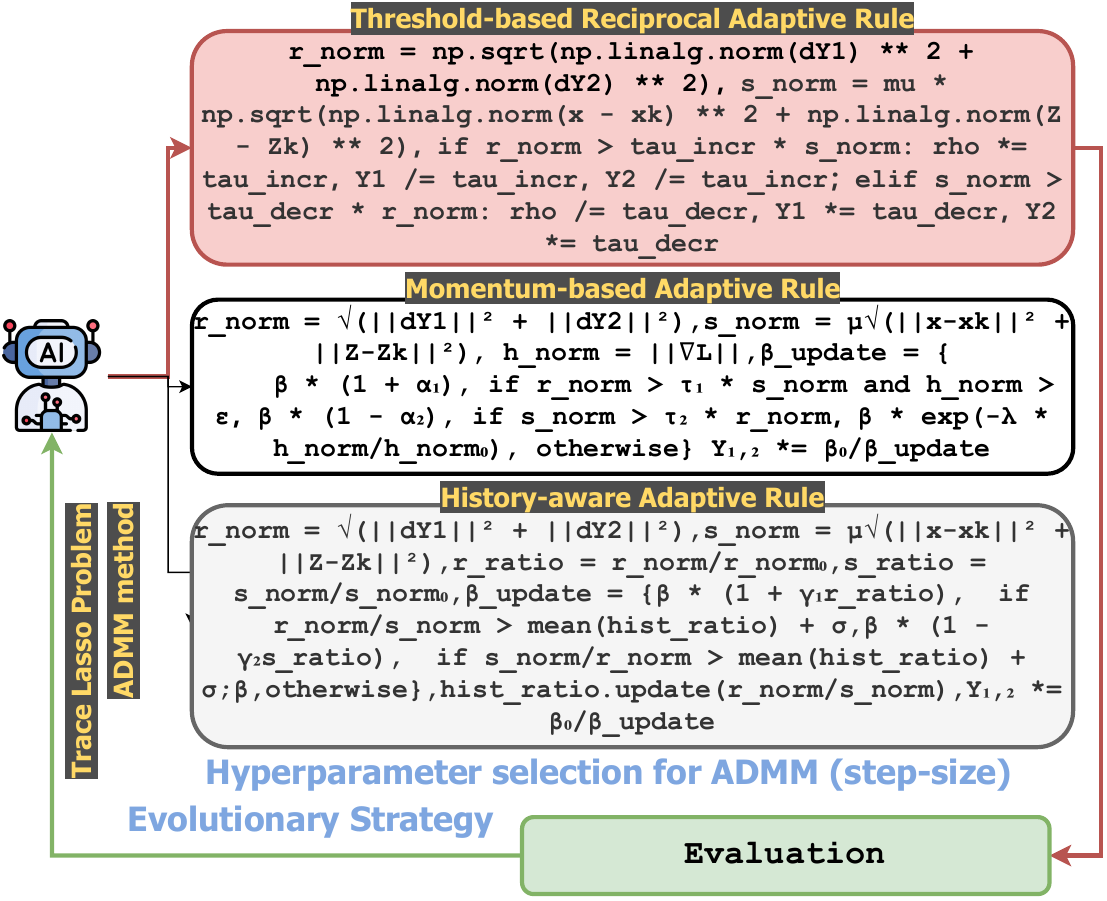}
    \caption{Evolutionary agentic workflow for ADMM.}
    \label{fig:admm}
\end{figure}

Recall the structured convex optimization problem \eqref{eq:convex}, while the iteration scheme of the classical ADMM method is denoted as (the $k$-th iteration),
\begin{equation}\label{eq:ADMM-scheme}
    \left\{\begin{array}{ll}
       \bx^{k} &= \arg\min_{\bx} {\cal{L}} (\bx,\bz^{k-1},\lambda^{k-1}) \\ 
       \bz^{k} &= \arg\min_{\bz} {\cal{L}} (\bx^k,\bz,\lambda^{k-1}) \\ 
       \lambda^k &= \lambda^{k-1} - \beta \left( A \bx^k + B \bz^k - c \right),
    \end{array}
    \right.
\end{equation}where
\begin{equation}\label{eq:Lagrangian}
\begin{array}{ll}
{\cal{L}} (\bx,\bz,\lambda) &= f(\bx) + g(\bz) - \lambda^{\rm T} \left( A \bx + B \bz - c \right) \\
& \qquad\qquad \qquad\qquad + \frac{\beta}{2}\left\| A \bx + B \bz - c \right\|^2_2,
\end{array}
\end{equation} where $\lambda$ denotes the dual variable (or Lagrange multiplier) concerning the linear constraint, and $\beta$ denotes the penalty parameter.
In general, the parameter $\beta$ is set to be a constant.

As a representative primal-dual method, the computational efficiency or convergence speed is closely related not only to primal variables $\bx$ and $\bz$ but the dual variable $\lambda$.
Typically, balancing the convergence speed of the primal variable and the dual variable is natural to guarantee better computational efficiency.
As in \cite{he2000alternating}, we define the primal $p_k$ and dual residual $d_k$ in the $k$-th iteration as follows:
\begin{equation}\label{eq:residual}
    \left\{ \begin{array}{ll}
         p_k &= \rho A^{\rm T} B \left( \bz^k - \bz^{k-1} \right), \\
         d_k &= A \bx^{k} + B \bz^{k} - c.\\  
    \end{array} \right.
\end{equation}

As highlighted in many previous works~\cite{he2000alternating,xu2017adaptive,xu2017admm},
the computational efficiency or convergence speed of ADMM is highly sensitive to the choice of the parameter $\beta$.
As an initial exploration that differs from the fixed constant parameter setting, \cite{he2000alternating} novelly proposed a self-adaptive scheme to adjust dynamically $\beta$ by balancing the primal residua $p_k$ and dual residual $d_k$.
In the $(k+1)$-th iteration, the parameter $\beta_{k+1}$ is determined by the relationship between the primal residua $p_k$ and dual residual $d_k$.
As a result, the parameter $\beta$ can be formalized as a function 
$$\beta_{k+1} = h \left( p_k, d_k \right).$$
For instance, \cite{he2000alternating} introduced an expert-designed nonlinear function, i.e.,
\begin{equation}\label{eq:self-adaptive}
    \beta_{k+1} = h \left( p_k, d_k \right) = \left\{
    \begin{array}{ll}
        \eta \beta_k, & \left\| d_k \right\|_2 > \mu ||p_k||_2; \\
        \beta_k / \eta, & ||p_k||_2 > \mu ||d_k||_2; \\
        \beta_k, & \text{otherwise},
    \end{array}
    \right.
\end{equation}where $\mu > 1$ and $\eta > 1$ are hyperparameters.
It has been proven in many application problems that the above self-adaptive parameter tuning scheme performs much better than the constant parameter setting.
However, we cannot conclude that the scheme \eqref{eq:self-adaptive} is the ``optimal" function $h$.
In other words, this scheme should be problem-driven, and general schemes designed by experts may not necessarily perform well on specific problems.
This provides us with greater possibilities to search for better strategies.


We can establish a novel algorithm innovation framework following the evolutionary agentic workflow, which uses large language models (LLMs) to explore self-adaptive strategies for the parameter automatically $\beta$, as shown in Figure~\ref{fig:admm}.
We can employ a pre-trained LLM for generating and searching the self-adaptive parameter scheme instead of designing the tuning scheme by experts as \eqref{eq:self-adaptive}.


To prove the possibility and efficiency of the proposed framework, we 
test on a Library of ADMM for sparse and low-rank
optimization, called LibADMM~\footnote{\url{https://github.com/canyilu/LibADMM}.}
All the problems in LibADMM can be calculated using the classical ADMM.

The new hyperparameter tuning criteria obtained through the agent workflow have led to the construction of an enhanced ADMM method, named NewADMM. We first demonstrate the NewADMM method's computational capabilities by training and testing on eight problems from LibADMM. As shown in Table \ref{tab: admm-general}, the results reveal that compared to the ADMM-expert approach, the NewADMM method requires significantly fewer iterations to achieve the same computational accuracy.

\begin{table}[!ht]
    \centering
    \resizebox{.7\linewidth}{!}{%
    \begin{tabular}{c|c|c|c|c}
    \hline
        Problem & L1 & Elasticnet & L1R & ElasticnetR \\ \hline
        ADMM-expert & 227 & 229 & 180 & 250  \\ \hline
        \textbf{\textcolor{purple}{NewADMM}} & \textbf{\textcolor{purple}{60}} & \textbf{\textcolor{purple}{12}} & \textbf{\textcolor{purple}{179}} & \textbf{\textcolor{purple}{24}}  \\ \hline
        \hline
        Problem & RPCA & LRMC & LRR & RMSC \\ \hline
        ADMM-expert & 94 & 83 & 198 & 116 \\ \hline
        \textbf{\textcolor{purple}{NewADMM}} & \textbf{\textcolor{purple}{20}} & \textbf{\textcolor{purple}{122}} & \textbf{\textcolor{purple}{3}} & \textbf{\textcolor{purple}{24}} \\ \hline
    \end{tabular}%
    }
    \label{tab: admm-general}
    \caption{Iterative number comparison on LibADMM applications to abtain the same tolerence.}
\end{table}

Further, in Table \ref{tab: admm-generalization}, we test the generalization capability of the new hyperparameter tuning criteria. 
Different from the results in Table \ref{tab: admm-general}, we specifically designed the hyperparameter tuning criteria by optimizing them on the L1 and Elasticnet problems within the LibADMM application suite, then directly transferred these criteria to other issues in LibADMM without further fine-tuning. 
The results demonstrate that when applying the criteria derived from the L1 and Elasticnet problems to other tasks, the same computational accuracy can still be achieved with significantly fewer iteration steps. 
This further validates the effectiveness of the agent workflow technology in developing robust hyperparameter tuning criteria.

\begin{table}[!ht]
    \centering
    \resizebox{.7\linewidth}{!}{%
    \begin{tabular}{c|c|c|c}
    \hline
        ~ & Problem & ADMM-expert & \textbf{\textcolor{purple}{NewADMM}} \\ \hline
        \multirow{6}*{L1} & Elasticnet & 235 & \textbf{\textcolor{purple}{45}} \\ \cline{2-4}
         & GroupL1 & 185 & \textbf{\textcolor{purple}{39}} \\ \cline{2-4}
         & L1R & 187 & \textbf{\textcolor{purple}{35}} \\ \cline{2-4}
         & ElasticnetR & 255 & \textbf{\textcolor{purple}{52}} \\ \cline{2-4}
         & TraceLasso & 92 & \textbf{\textcolor{purple}{43}} \\ \cline{2-4}
         & GroupL1R & 187 & \textbf{\textcolor{purple}{40}} \\ \hline
        \hline
        \multirow{6}*{Elasticnet} & L1 & 237 & \textbf{\textcolor{purple}{36}} \\ \cline{2-4}
         & GroupL1 & 185 & \textbf{\textcolor{purple}{34}} \\ \cline{2-4}
         & L1R & 187 & \textbf{\textcolor{purple}{30}} \\ \cline{2-4}
         & ElasticnetR & 255 & \textbf{\textcolor{purple}{154}} \\ \cline{2-4}
         & TraceLasso & 92 & \textbf{\textcolor{purple}{52}} \\ \cline{2-4}
         & GroupL1R & 187 & \textbf{\textcolor{purple}{34}} \\ \hline
    \end{tabular}%
    }
    \label{tab: admm-generalization}
    \caption{Performance of generalization.}
\end{table}

\section{Limitations}

While the proposed evolutionary agentic workflow shows promise in revolutionizing optimization practice, there remain two key challenges with promising solutions on the horizon. 
First, foundation agents lack built-in mechanisms for theoretical verification. 
In optimization research, theoretical guarantees such as convergence properties and optimality conditions are fundamental requirements, not optional features. 
When agents propose novel optimization strategies through evolution, their mathematical soundness remains unverified, potentially limiting academic adoption.

Second, the high inference cost of FMs presents a significant barrier to industrial adoption. 
While the proposed framework aims to bridge the academic-industrial gap by lowering the expertise barrier, the computational overhead of repeated agent interactions may make it prohibitively expensive for practical deployment, especially in industrial settings where cost-efficiency is crucial.

Third, scalability to large-scale industrial optimization problems represents a significant practical challenge. 
While the framework offers several advantages, the computational demands of evaluating solutions for huge problem instances could become prohibitive without appropriate infrastructure support.

For theoretical verification, integrating code FMs~\citep{roziere2023code,jiang2024survey} with automated proof assistants~\citep{de2015lean,moura2021lean,song2024towards} could enable automated verification of mathematical properties.
Efforts in formalizing block-partitioned methods, KKT conditions, and subdifferential-based convergence~\citep{li2025formalization1,li2025formalization2} suggest automatically verifying ADMM and other iterative schemes. 
Future work involves integrating Lean4 directly into the proposed agentic workflow, enabling proofs and generation of optimization code in a single loop.

For inference cost, advances in model quantization~\citep{lin2024awq,egashira2024exploiting}, distillation~\citep{xu2024survey}, and small FMs~\citep{schick2020s,van2024survey} offer promising solutions for cost-effective industrial deployment.

Regarding scalability challenges, several promising approaches could address these concerns. 
The proposed workflow is inherently heuristic-focused, optimizing approaches not restricted by problem size. 
The primary bottleneck emerges when validating solutions for huge instances, which could be addressed through high-performance simulation or distributed HPC frameworks. 
Furthermore, bottom-up operator optimization—automating kernel generation and leveraging specialized accelerators (e.g., GPU or HPC cloud resources)—could drastically expedite the evaluation phase for large-scale tasks. 
This approach parallels ongoing research efforts where agentic workflows have generated optimized CUDA kernels to accelerate large-model inference, and similar speedups could benefit large-scale industrial optimization~\citep{lange2025ai,chen2025automating}.

Addressing these limitations would move us toward a unified framework that maintains theoretical rigor while enabling practical deployment, ultimately making advanced optimization techniques more accessible to practitioners.

\section{Alternative Views}

Key components of our framework face substantive challenges from existing research. 
Studies on LLMs' mathematical capabilities suggest fundamental limitations in complex reasoning. 
Recent work has demonstrated that these models often fail at multi-step mathematical derivations and struggle to maintain logical consistency in complex optimization problems~\citep{mirzadeh2024gsm,yadkori2024believe,kambhampati2024can,valmeekam2022large,valmeekam2024llms,boix-adsera2024when,jiang2024peek}, raising concerns about their reliability in algorithm design and parameter tuning.

The evolutionary aspect of our approach faces more fundamental criticisms. 
Traditional evolutionary computation methods suffer from inherent limitations in rigorous performance analysis~\citep{kudela2022critical,kudela2023evolutionary}. 
Many newly proposed methods show superior performance primarily due to implicit biases in benchmark functions, such as having optima at the center of the feasible set.
Furthermore, when tested on more challenging scenarios like real-world applications, many evolutionary methods perform no better than classical approaches or even random search~\citep{wolpert1997no,araujo2007evolutionary,woldesenbet2009dynamic,yampolskiy2018we,velasco2024literature}.

Recent advances suggest these challenges can be effectively addressed. 
LLM trained via large-scale reinforcement learning has shown enhanced mathematical reasoning capabilities~\citep{openai2024o1,guo2025deepseek}. At the same time, the agentic workflow compensates for evolutionary methods' limitations by incorporating domain knowledge and efficient solution representations through FMs. 
This synergy creates a more principled optimization framework that leverages improved reasoning capabilities and structured search strategies.

\section{Related Work}

The automation of optimization processes has been explored from various perspectives. 
Automated Machine Learning (AutoML) and Learning to Optimize (L2O) are two specific methodologies focusing on ML tasks and optimization algorithms. 
Meta-optimization and Automatic Algorithm Configuration (AAC) offer more general frameworks, addressing the higher-level challenge of optimizing optimization processes. 
While these approaches have succeeded significantly in their respective domains, they differ from evolutionary agentic workflow in scope and methodology.

\subsection{Automated Machine Learning}
AutoML represents a significant step towards automating the ML pipeline~\citep{hutter2019automated,yao2018taking}. 
It primarily focuses on automating model selection, feature engineering, and hyperparameter optimization for machine learning tasks~\citep{feurer2019hyperparameter}, employing sophisticated search strategies such as Bayesian optimization~\citep{snoek2012practical}, evolutionary algorithms~\citep{real2019regularized}.
even agentic workflow~\citep{zhang2023automl,hong2024data,trirat2024automl,chi2024sela,gu2024large,hoseini2024enhancing}.
However, AutoML faces unique challenges when applied to general optimization problems. 
While machine learning tasks often share common structural patterns that can be effectively explored through predefined search spaces~\citep{liu2018darts,pham2018efficient,elsken2019neural}, optimization problems typically require more sophisticated mathematical formulations and theoretical guarantees~\citep{wolpert1997no,boyd2004convex}. 
Evolutionary agentic workflow potentially extends beyond AutoML's fixed search space approach by potentially leveraging FMs to synthesize new optimization strategies and problem formulations.

\subsection{Learning to Optimize}

L2O attempts to automate optimization processes through end-to-end neural architectures~\citep{andrychowicz2016learning,chen2022learning}. 
This approach treats the optimization algorithm as a learnable component, typically using recurrent neural networks to generate update rules that can potentially outperform traditional optimization methods, showing particular strength in learning highly efficient optimization strategies for specific problem classes~\citep{li2017learning}.
While L2O focuses primarily on learning the optimization algorithm with limited generalization capability~\citep{lv2017learning}, evolutionary agentic workflow addresses the entire optimization pipeline.
Additionally, agentic workflow can integrate existing optimization theory and expert insights through FMs, overcoming L2O's limitation in incorporating domain knowledge~\citep{metz2019understanding}.

\subsection{Meta-Optimization}

Meta-optimization, initially developed in ML~\citep{finn2017model,devlin2017neural}, has expanded to general optimization problems~\citep{franceschi2018bilevel}. 
This paradigm focuses on optimizing the optimization process itself, achieving significant success in few-shot learning and transfer learning~\citep{hospedales2021meta}, while also showing promise in adapting solver configurations and learning initialization strategies~\citep{ma2024toward}.
The agentic framework distinguishes itself from meta-optimization in its approach to problem-solving. 
While meta-optimization typically requires a predefined family of optimization problems and focuses on learning transferable optimization strategies within this family, the agentic framework aims to understand and formulate optimization problems from natural language descriptions, enabling the generation of novel solution approaches rather than just adapting existing ones.

\subsection{Automatic Algorithm Configuration}

Automatic algorithm configuration (AAC) has evolved from its origins in operations research~\citep{hutter2011sequential,lopez2016irace} to a widely adopted methodology across multiple domains, including simulation software~\citep{bartz2020benchmarking}, machine learning pipelines~\citep{feurer2022auto}, and various algorithmic frameworks~\citep{hutter2017configurable,gleixner2021miplib}. 
Using sophisticated search strategies like iterated racing and model-based search, AAC has successfully tuned complex systems with numerous interacting parameters.
However, AAC's scope is fundamentally narrower than our proposed framework. 
While AAC excels at finding optimal configurations within a predefined parameter space, it cannot reformulate problems or synthesize new algorithmic components.
The agentic workflow extends beyond parameter tuning to enable innovation across the optimization pipeline, from problem formulation to algorithm design.

\section{Conclusion}

This paper advocates a paradigm shift in optimization problem solving through evolutionary agentic workflows. 
While existing research has demonstrated foundation agents' effectiveness in individual optimization components, our framework proposes their evolutionary integration into a comprehensive system. 
Although we acknowledge current verification and computational overhead limitations, recent advances suggest these challenges can be systematically addressed. 
Our framework combines FMs' proven capabilities with evolutionary methods' structured exploration, representing a crucial step in bridging the gap between academic innovation and industrial practice.


\bibliography{main}

\begin{thebibliography}{180}
\providecommand{\natexlab}[1]{#1}
\providecommand{\url}[1]{\texttt{#1}}
\expandafter\ifx\csname urlstyle\endcsname\relax
  \providecommand{\doi}[1]{doi: #1}\else
  \providecommand{\doi}{doi: \begingroup \urlstyle{rm}\Url}\fi

\bibitem[Achterberg \& Wunderling(2013)Achterberg and Wunderling]{achterberg2013mixed}
Achterberg, T. and Wunderling, R.
\newblock Mixed integer programming: Analyzing 12 years of progress.
\newblock In \emph{Facets of Combinatorial Optimization: Festschrift for Martin Gr{\"o}tschel}, pp.\  449--481. Springer, 2013.

\bibitem[AhmadiTeshnizi et~al.(2023)AhmadiTeshnizi, Gao, and Udell]{ahmaditeshnizi2023optimus}
AhmadiTeshnizi, A., Gao, W., and Udell, M.
\newblock Optimus: {O}ptimization modeling using mip solvers and large language models.
\newblock \emph{arXiv:2310.06116}, 2023.

\bibitem[AhmadiTeshnizi et~al.(2024)AhmadiTeshnizi, Gao, and Udell]{ahmaditeshnizioptimus}
AhmadiTeshnizi, A., Gao, W., and Udell, M.
\newblock Optimus: Scalable optimization modeling with (mi) lp solvers and large language models.
\newblock In \emph{ICML}, 2024.

\bibitem[Ahmed \& Choudhury(2024)Ahmed and Choudhury]{ahmed2024lm4opt}
Ahmed, T. and Choudhury, S.
\newblock {LM4OPT}: {U}nveiling the potential of large language models in formulating mathematical optimization problems.
\newblock \emph{INFOR: Information Systems and Operational Research}, 62\penalty0 (4):\penalty0 559--572, 2024.

\bibitem[Akiba et~al.(2019)Akiba, Sano, Yanase, Ohta, and Koyama]{akiba2019optuna}
Akiba, T., Sano, S., Yanase, T., Ohta, T., and Koyama, M.
\newblock {O}ptuna: A next-generation hyperparameter optimization framework.
\newblock In \emph{KDD}, 2019.

\bibitem[Almonacid(2023)]{almonacid2023towards}
Almonacid, B.
\newblock Towards an automatic optimisation model generator assisted with generative pre-trained transformer.
\newblock \emph{arXiv:2305.05811}, 2023.

\bibitem[Amarasinghe et~al.(2023)Amarasinghe, Nguyen, Sun, and Alahakoon]{amarasinghe2023ai}
Amarasinghe, P.~T., Nguyen, S., Sun, Y., and Alahakoon, D.
\newblock {AI-Copilot} for business optimisation: {A} framework and a case study in production scheduling.
\newblock \emph{arXiv:2309.13218}, 2023.

\bibitem[Andrychowicz et~al.(2016)Andrychowicz, Denil, Gomez, Hoffman, Pfau, Schaul, Shillingford, and De~Freitas]{andrychowicz2016learning}
Andrychowicz, M., Denil, M., Gomez, S., Hoffman, M.~W., Pfau, D., Schaul, T., Shillingford, B., and De~Freitas, N.
\newblock Learning to learn by gradient descent by gradient descent.
\newblock \emph{NeurIPS}, 2016.

\bibitem[Araujo(2007)]{araujo2007evolutionary}
Araujo, L.
\newblock How evolutionary algorithms are applied to statistical natural language processing.
\newblock \emph{Artificial Intelligence Review}, 28:\penalty0 275--303, 2007.

\bibitem[Bai et~al.(2022)Bai, Kadavath, Kundu, Askell, Kernion, Jones, Chen, Goldie, Mirhoseini, McKinnon, et~al.]{bai2022constitutional}
Bai, Y., Kadavath, S., Kundu, S., Askell, A., Kernion, J., Jones, A., Chen, A., Goldie, A., Mirhoseini, A., McKinnon, C., et~al.
\newblock Constitutional ai: Harmlessness from ai feedback.
\newblock \emph{arXiv:2212.08073}, 2022.

\bibitem[Bartz-Beielstein et~al.(2020)Bartz-Beielstein, Doerr, Berg, Bossek, Chandrasekaran, Eftimov, Fischbach, Kerschke, La~Cava, Lopez-Ibanez, et~al.]{bartz2020benchmarking}
Bartz-Beielstein, T., Doerr, C., Berg, D. v.~d., Bossek, J., Chandrasekaran, S., Eftimov, T., Fischbach, A., Kerschke, P., La~Cava, W., Lopez-Ibanez, M., et~al.
\newblock Benchmarking in optimization: Best practice and open issues.
\newblock \emph{arXiv:2007.03488}, 2020.

\bibitem[Bays(1977)]{bays1977comparison}
Bays, C.
\newblock A comparison of next-fit, first-fit, and best-fit.
\newblock \emph{Communications of the ACM}, 20\penalty0 (3):\penalty0 191--192, 1977.

\bibitem[Bazaraa \& Sherali(1981)Bazaraa and Sherali]{bazaraa1981choice}
Bazaraa, M.~S. and Sherali, H.~D.
\newblock On the choice of step-size in subgradient optimization.
\newblock \emph{European Journal of Operational Research}, 7\penalty0 (4):\penalty0 380--388, 1981.

\bibitem[Bernal-Agust{\'\i}n \& Dufo-Lopez(2009)Bernal-Agust{\'\i}n and Dufo-Lopez]{bernal2009simulation}
Bernal-Agust{\'\i}n, J.~L. and Dufo-Lopez, R.
\newblock Simulation and optimization of stand-alone hybrid renewable energy systems.
\newblock \emph{Renewable and sustainable energy reviews}, 13\penalty0 (8):\penalty0 2111--2118, 2009.

\bibitem[Bixby(2002)]{bixby2002solving}
Bixby, R.~E.
\newblock Solving real-world linear programs: {A} decade and more of progress.
\newblock \emph{Operations Research}, 50\penalty0 (1):\penalty0 3--15, 2002.

\bibitem[Bohanec \& Bratko(1994)Bohanec and Bratko]{bohanec1994trading}
Bohanec, M. and Bratko, I.
\newblock Trading accuracy for simplicity in decision trees.
\newblock \emph{Machine Learning}, 15:\penalty0 223--250, 1994.

\bibitem[Boix-Adser{\`a} et~al.(2024)Boix-Adser{\`a}, Saremi, Abbe, Bengio, Littwin, and Susskind]{boix-adsera2024when}
Boix-Adser{\`a}, E., Saremi, O., Abbe, E., Bengio, S., Littwin, E., and Susskind, J.~M.
\newblock When can transformers reason with abstract symbols?
\newblock In \emph{ICLR}, 2024.

\bibitem[Boyd \& Mattingley(2007)Boyd and Mattingley]{boyd2007branch}
Boyd, S. and Mattingley, J.
\newblock Branch and bound methods.
\newblock \emph{Notes for EE364b, Stanford University}, 2006:\penalty0 07, 2007.

\bibitem[Boyd \& Vandenberghe(2004)Boyd and Vandenberghe]{boyd2004convex}
Boyd, S. and Vandenberghe, L.
\newblock Convex optimization.
\newblock \emph{Cambridge University Press}, 2004.

\bibitem[Boyd et~al.(2011)Boyd, Parikh, Chu, Peleato, Eckstein, et~al.]{boyd2011distributed}
Boyd, S., Parikh, N., Chu, E., Peleato, B., Eckstein, J., et~al.
\newblock Distributed optimization and statistical learning via the alternating direction method of multipliers.
\newblock \emph{Foundations and Trends{\textregistered} in Machine learning}, 3\penalty0 (1):\penalty0 1--122, 2011.

\bibitem[Chan et~al.(2024)Chan, Chen, Su, Yu, Xue, Zhang, Fu, and Liu]{chan2024chateval}
Chan, C.-M., Chen, W., Su, Y., Yu, J., Xue, W., Zhang, S., Fu, J., and Liu, Z.
\newblock Chateval: {T}owards better {LLM}-based evaluators through multi-agent debate.
\newblock In \emph{ICLR}, 2024.

\bibitem[Chang et~al.(2020)Chang, Hong, Wai, Zhang, and Lu]{chang2020distributed}
Chang, T.-H., Hong, M., Wai, H.-T., Zhang, X., and Lu, S.
\newblock Distributed learning in the nonconvex world: From batch data to streaming and beyond.
\newblock \emph{IEEE Signal Processing Magazine}, 37\penalty0 (3):\penalty0 26--38, 2020.

\bibitem[Chen et~al.(2023{\natexlab{a}})Chen, Shu, Shareghi, Collier, Narasimhan, and Yao]{chen2023fireact}
Chen, B., Shu, C., Shareghi, E., Collier, N., Narasimhan, K., and Yao, S.
\newblock Fireact: Toward language agent fine-tuning.
\newblock \emph{arXiv preprint arXiv:2310.05915}, 2023{\natexlab{a}}.

\bibitem[Chen et~al.(2024)Chen, Constante-Flores, and Li]{chen2024diagnosing}
Chen, H., Constante-Flores, G.~E., and Li, C.
\newblock Diagnosing infeasible optimization problems using large language models.
\newblock \emph{INFOR: Information Systems and Operational Research}, 62\penalty0 (4):\penalty0 573--587, 2024.

\bibitem[Chen et~al.(2022)Chen, Chen, Chen, Heaton, Liu, Wang, and Yin]{chen2022learning}
Chen, T., Chen, X., Chen, W., Heaton, H., Liu, J., Wang, Z., and Yin, W.
\newblock Learning to optimize: {A} primer and a benchmark.
\newblock \emph{Journal of Machine Learning Research}, 23\penalty0 (189):\penalty0 1--59, 2022.

\bibitem[Chen et~al.(2025)Chen, Xu, and Devleker]{chen2025automating}
Chen, T., Xu, B., and Devleker, K.
\newblock Automating gpu kernel generation with deepseek-r1 and inference time scaling, February 2025.

\bibitem[Chen et~al.(2023{\natexlab{b}})Chen, Deng, Shen, and Yu]{chen2023mind}
Chen, X., Deng, T., Shen, Z.-J.~M., and Yu, Y.
\newblock Mind the gap between research and practice in operations management.
\newblock \emph{IISE Transactions}, 55\penalty0 (1):\penalty0 32--42, 2023{\natexlab{b}}.

\bibitem[Cheng et~al.(2024)Cheng, Zhang, Zhang, Meng, Hong, Li, Wang, Wang, Yin, Zhao, et~al.]{cheng2024exploring}
Cheng, Y., Zhang, C., Zhang, Z., Meng, X., Hong, S., Li, W., Wang, Z., Wang, Z., Yin, F., Zhao, J., et~al.
\newblock Exploring large language model based intelligent agents: {D}efinitions, methods, and prospects.
\newblock \emph{arXiv:2401.03428}, 2024.

\bibitem[Chi et~al.(2024)Chi, Lin, Hong, Pan, Fei, Mei, Liu, Pang, Kwok, Zhang, et~al.]{chi2024sela}
Chi, Y., Lin, Y., Hong, S., Pan, D., Fei, Y., Mei, G., Liu, B., Pang, T., Kwok, J., Zhang, C., et~al.
\newblock {SELA}: {T}ree-search enhanced llm agents for automated machine learning.
\newblock \emph{arXiv:2410.17238}, 2024.

\bibitem[Chong \& {\.Z}ak(2013)Chong and {\.Z}ak]{chong2013introduction}
Chong, E.~K. and {\.Z}ak, S.~H.
\newblock \emph{An Introduction to Optimization}, volume~75.
\newblock John Wiley \& Sons, 2013.

\bibitem[Chu et~al.(2023)Chu, Chen, Chen, Yu, He, Wang, Peng, Liu, Qin, and Liu]{chu2023survey}
Chu, Z., Chen, J., Chen, Q., Yu, W., He, T., Wang, H., Peng, W., Liu, M., Qin, B., and Liu, T.
\newblock A survey of chain of thought reasoning: {A}dvances, frontiers and future.
\newblock \emph{arXiv:2309.15402}, 2023.

\bibitem[C{\^o}t{\'e} et~al.(2021)C{\^o}t{\'e}, Haouari, and Iori]{cote2021combinatorial}
C{\^o}t{\'e}, J.-F., Haouari, M., and Iori, M.
\newblock Combinatorial benders decomposition for the two-dimensional bin packing problem.
\newblock \emph{INFORMS Journal on Computing}, 33\penalty0 (3):\penalty0 963--978, 2021.

\bibitem[Dantzig(2002)]{dantzig2002linear}
Dantzig, G.~B.
\newblock Linear programming.
\newblock \emph{Operations Research}, 50\penalty0 (1):\penalty0 42--47, 2002.

\bibitem[De~Moura et~al.(2015)De~Moura, Kong, Avigad, Van~Doorn, and von Raumer]{de2015lean}
De~Moura, L., Kong, S., Avigad, J., Van~Doorn, F., and von Raumer, J.
\newblock The {L}ean theorem prover (system description).
\newblock In \emph{CADE}, 2015.

\bibitem[Devlin et~al.(2017)Devlin, Bunel, Singh, Hausknecht, and Kohli]{devlin2017neural}
Devlin, J., Bunel, R.~R., Singh, R., Hausknecht, M., and Kohli, P.
\newblock Neural program meta-induction.
\newblock \emph{NeurIPS}, 2017.

\bibitem[Dimopoulos \& Zalzala(2000)Dimopoulos and Zalzala]{dimopoulos2000recent}
Dimopoulos, C. and Zalzala, A.~M.
\newblock Recent developments in evolutionary computation for manufacturing optimization: problems, solutions, and comparisons.
\newblock \emph{IEEE Transactions on Evolutionary Computation}, 4\penalty0 (2):\penalty0 93--113, 2000.

\bibitem[Dincer et~al.(2017)Dincer, Rosen, and Ahmadi]{dincer2017optimization}
Dincer, I., Rosen, M.~A., and Ahmadi, P.
\newblock \emph{Optimization of Energy Systems}.
\newblock John Wiley \& Sons, 2017.

\bibitem[Egashira et~al.(2024)Egashira, Vero, Staab, He, and Vechev]{egashira2024exploiting}
Egashira, K., Vero, M., Staab, R., He, J., and Vechev, M.
\newblock Exploiting {LLM} quantization.
\newblock \emph{arXiv:2405.18137}, 2024.

\bibitem[Elsken et~al.(2019)Elsken, Metzen, and Hutter]{elsken2019neural}
Elsken, T., Metzen, J.~H., and Hutter, F.
\newblock Neural architecture search: A survey.
\newblock \emph{Journal of Machine Learning Research}, 20\penalty0 (55):\penalty0 1--21, 2019.

\bibitem[Fan et~al.(2024)Fan, Ghaddar, Wang, Xing, Zhang, and Zhou]{fan2024artificial}
Fan, Z., Ghaddar, B., Wang, X., Xing, L., Zhang, Y., and Zhou, Z.
\newblock Artificial intelligence for operations research: {R}evolutionizing the operations research process.
\newblock \emph{arXiv:2401.03244}, 2024.

\bibitem[Feurer \& Hutter(2019)Feurer and Hutter]{feurer2019hyperparameter}
Feurer, M. and Hutter, F.
\newblock Hyperparameter optimization.
\newblock \emph{Automated machine learning: Methods, systems, challenges}, pp.\  3--33, 2019.

\bibitem[Feurer et~al.(2022)Feurer, Eggensperger, Falkner, Lindauer, and Hutter]{feurer2022auto}
Feurer, M., Eggensperger, K., Falkner, S., Lindauer, M., and Hutter, F.
\newblock Auto-sklearn 2.0: {H}ands-free automl via meta-learning.
\newblock \emph{Journal of Machine Learning Research}, 23\penalty0 (261):\penalty0 1--61, 2022.

\bibitem[Finn et~al.(2017)Finn, Abbeel, and Levine]{finn2017model}
Finn, C., Abbeel, P., and Levine, S.
\newblock Model-agnostic meta-learning for fast adaptation of deep networks.
\newblock In \emph{ICML}, 2017.

\bibitem[Franceschi et~al.(2018)Franceschi, Frasconi, Salzo, Grazzi, and Pontil]{franceschi2018bilevel}
Franceschi, L., Frasconi, P., Salzo, S., Grazzi, R., and Pontil, M.
\newblock Bilevel programming for hyperparameter optimization and meta-learning.
\newblock In \emph{ICML}, 2018.

\bibitem[Gangwar \& Kani(2023)Gangwar and Kani]{gangwar2023highlighting}
Gangwar, N. and Kani, N.
\newblock Highlighting named entities in input for auto-formulation of optimization problems.
\newblock In \emph{International Conference on Intelligent Computer Mathematics}, 2023.

\bibitem[Ge et~al.(2023)Ge, Hua, Mei, Tan, Xu, Li, Zhang, et~al.]{ge2023openagi}
Ge, Y., Hua, W., Mei, K., Tan, J., Xu, S., Li, Z., Zhang, Y., et~al.
\newblock Open{AGI}: When {LLM} meets domain experts.
\newblock \emph{NeurIPS}, 2023.

\bibitem[Giovanelli et~al.(2024)Giovanelli, Tornede, Tornede, and Lindauer]{giovanelli2024interactive}
Giovanelli, J., Tornede, A., Tornede, T., and Lindauer, M.
\newblock Interactive hyperparameter optimization in multi-objective problems via preference learning.
\newblock In \emph{AAAI}, 2024.

\bibitem[Gleixner et~al.(2021)Gleixner, Hendel, Gamrath, Achterberg, Bastubbe, Berthold, Christophel, Jarck, Koch, Linderoth, et~al.]{gleixner2021miplib}
Gleixner, A., Hendel, G., Gamrath, G., Achterberg, T., Bastubbe, M., Berthold, T., Christophel, P., Jarck, K., Koch, T., Linderoth, J., et~al.
\newblock {MIPLIB} 2017: {D}ata-driven compilation of the 6th mixed-integer programming library.
\newblock \emph{Mathematical Programming Computation}, 13\penalty0 (3):\penalty0 443--490, 2021.

\bibitem[Greenhalgh \& Marshall(2000)Greenhalgh and Marshall]{greenhalgh2000convergence}
Greenhalgh, D. and Marshall, S.
\newblock Convergence criteria for genetic algorithms.
\newblock \emph{SIAM Journal on Computing}, 30\penalty0 (1):\penalty0 269--282, 2000.

\bibitem[Gu et~al.(2024)Gu, You, Cao, and Yu]{gu2024large}
Gu, Y., You, H., Cao, J., and Yu, M.
\newblock Large language models for constructing and optimizing machine learning workflows: {A} survey.
\newblock \emph{arXiv:2411.10478}, 2024.

\bibitem[Guo et~al.(2025)Guo, Yang, Zhang, Song, Zhang, Xu, Zhu, Ma, Wang, Bi, et~al.]{guo2025deepseek}
Guo, D., Yang, D., Zhang, H., Song, J., Zhang, R., Xu, R., Zhu, Q., Ma, S., Wang, P., Bi, X., et~al.
\newblock Deepseek-r1: Incentivizing reasoning capability in llms via reinforcement learning.
\newblock \emph{arXiv:2501.12948}, 2025.

\bibitem[Hao et~al.(2023)Hao, Gu, Ma, Hong, Wang, Wang, and Hu]{hao2023reasoning}
Hao, S., Gu, Y., Ma, H., Hong, J., Wang, Z., Wang, D., and Hu, Z.
\newblock Reasoning with language model is planning with world model.
\newblock In \emph{EMNLP}, 2023.

\bibitem[He et~al.(2000)He, Yang, and Wang]{he2000alternating}
He, B., Yang, H., and Wang, S.
\newblock Alternating direction method with self-adaptive penalty parameters for monotone variational inequalities.
\newblock \emph{Journal of Optimization Theory and applications}, 106:\penalty0 337--356, 2000.

\bibitem[He et~al.(2016)He, Ma, and Yuan]{he2016convergence}
He, B., Ma, F., and Yuan, X.
\newblock Convergence study on the symmetric version of {ADMM} with larger step sizes.
\newblock \emph{SIAM Journal on Imaging Sciences}, 9\penalty0 (3):\penalty0 1467--1501, 2016.

\bibitem[Hong et~al.(2015)Hong, Razaviyayn, Luo, and Pang]{hong2015unified}
Hong, M., Razaviyayn, M., Luo, Z.-Q., and Pang, J.-S.
\newblock A unified algorithmic framework for block-structured optimization involving big data: With applications in machine learning and signal processing.
\newblock \emph{IEEE Signal Processing Magazine}, 33\penalty0 (1):\penalty0 57--77, 2015.

\bibitem[Hong et~al.(2024)Hong, Lin, Liu, Liu, Wu, Zhang, Wei, Li, Chen, Zhang, et~al.]{hong2024data}
Hong, S., Lin, Y., Liu, B., Liu, B., Wu, B., Zhang, C., Wei, C., Li, D., Chen, J., Zhang, J., et~al.
\newblock Data interpreter: {A}n {LLM} agent for data science.
\newblock \emph{arXiv:2402.18679}, 2024.

\bibitem[Hoseini et~al.(2024)Hoseini, Ibbels, and Quix]{hoseini2024enhancing}
Hoseini, S., Ibbels, M., and Quix, C.
\newblock Enhancing machine learning capabilities in data lakes with {AutoML} and {LLMs}.
\newblock In \emph{ADBIS}, 2024.

\bibitem[Hospedales et~al.(2021)Hospedales, Antoniou, Micaelli, and Storkey]{hospedales2021meta}
Hospedales, T., Antoniou, A., Micaelli, P., and Storkey, A.
\newblock Meta-learning in neural networks: {A} survey.
\newblock \emph{IEEE Transactions on Pattern Analysis and Machine Intelligence}, 44\penalty0 (9):\penalty0 5149--5169, 2021.

\bibitem[Huang et~al.(2024{\natexlab{a}})Huang, Yang, Qi, and Wang]{huang2024large}
Huang, S., Yang, K., Qi, S., and Wang, R.
\newblock When large language model meets optimization.
\newblock \emph{arXiv:2405.10098}, 2024{\natexlab{a}}.

\bibitem[Huang et~al.(2022)Huang, Xia, Xiao, Chan, Liang, Florence, Zeng, Tompson, Mordatch, Chebotar, Sermanet, Jackson, Brown, Luu, Levine, Hausman, and brian ichter]{huang2022inner}
Huang, W., Xia, F., Xiao, T., Chan, H., Liang, J., Florence, P., Zeng, A., Tompson, J., Mordatch, I., Chebotar, Y., Sermanet, P., Jackson, T., Brown, N., Luu, L., Levine, S., Hausman, K., and brian ichter.
\newblock Inner monologue: {E}mbodied reasoning through planning with language models.
\newblock In \emph{CoRL}, 2022.

\bibitem[Huang et~al.(2024{\natexlab{b}})Huang, Shen, Hu, Gao, and Wang]{huang2024mamo}
Huang, X., Shen, Q., Hu, Y., Gao, A., and Wang, B.
\newblock {MAMO}: a mathematical modeling benchmark with solvers.
\newblock \emph{arXiv:2405.13144}, 2024{\natexlab{b}}.

\bibitem[Hutter et~al.(2011)Hutter, Hoos, and Leyton-Brown]{hutter2011sequential}
Hutter, F., Hoos, H.~H., and Leyton-Brown, K.
\newblock Sequential model-based optimization for general algorithm configuration.
\newblock In \emph{LION}, 2011.

\bibitem[Hutter et~al.(2017)Hutter, Lindauer, Balint, Bayless, Hoos, and Leyton-Brown]{hutter2017configurable}
Hutter, F., Lindauer, M., Balint, A., Bayless, S., Hoos, H., and Leyton-Brown, K.
\newblock The configurable {SAT} solver challenge ({CSSC}).
\newblock \emph{Artificial Intelligence}, 243:\penalty0 1--25, 2017.

\bibitem[Hutter et~al.(2019)Hutter, Kotthoff, and Vanschoren]{hutter2019automated}
Hutter, F., Kotthoff, L., and Vanschoren, J.
\newblock \emph{Automated machine learning: methods, systems, challenges}.
\newblock Springer Nature, 2019.

\bibitem[Ibanez et~al.(2018)Ibanez, Clark, Huckman, and Staats]{ibanez2018discretionary}
Ibanez, M.~R., Clark, J.~R., Huckman, R.~S., and Staats, B.~R.
\newblock Discretionary task ordering: Queue management in radiological services.
\newblock \emph{Management Science}, 64\penalty0 (9):\penalty0 4389--4407, 2018.

\bibitem[Izquierdo et~al.(1999)Izquierdo, Medina, Vianna, Izquierdo, and Barros]{izquierdo1999separate}
Izquierdo, I., Medina, J.~H., Vianna, M.~R., Izquierdo, L.~A., and Barros, D.~M.
\newblock Separate mechanisms for short-and long-term memory.
\newblock \emph{Behavioural brain research}, 103\penalty0 (1):\penalty0 1--11, 1999.

\bibitem[Jiang et~al.(2024{\natexlab{a}})Jiang, Xie, Hao, Wang, Mallick, Su, Taylor, and Roth]{jiang2024peek}
Jiang, B., Xie, Y., Hao, Z., Wang, X., Mallick, T., Su, W.~J., Taylor, C.~J., and Roth, D.
\newblock A peek into token bias: {L}arge language models are not yet genuine reasoners.
\newblock \emph{arXiv:2406.11050}, 2024{\natexlab{a}}.

\bibitem[Jiang et~al.(2024{\natexlab{b}})Jiang, Wang, Shen, Kim, and Kim]{jiang2024survey}
Jiang, J., Wang, F., Shen, J., Kim, S., and Kim, S.
\newblock A survey on large language models for code generation.
\newblock \emph{arXiv:2406.00515}, 2024{\natexlab{b}}.

\bibitem[Jiang et~al.(2021)Jiang, Cao, and Zhang]{jiang2021learning}
Jiang, Y., Cao, Z., and Zhang, J.
\newblock Learning to solve 3-{D} bin packing problem via deep reinforcement learning and constraint programming.
\newblock \emph{IEEE Transactions on Cybernetics}, 53\penalty0 (5):\penalty0 2864--2875, 2021.

\bibitem[Kambhampati(2024)]{kambhampati2024can}
Kambhampati, S.
\newblock Can large language models reason and plan?
\newblock \emph{Annals of the New York Academy of Sciences}, 1534\penalty0 (1):\penalty0 15--18, 2024.

\bibitem[Khot et~al.(2023)Khot, Trivedi, Finlayson, Fu, Richardson, Clark, and Sabharwal]{khot2023decomposed}
Khot, T., Trivedi, H., Finlayson, M., Fu, Y., Richardson, K., Clark, P., and Sabharwal, A.
\newblock Decomposed prompting: {A} modular approach for solving complex tasks.
\newblock In \emph{ICLR}, 2023.

\bibitem[Kudela(2022)]{kudela2022critical}
Kudela, J.
\newblock A critical problem in benchmarking and analysis of evolutionary computation methods.
\newblock \emph{Nature Machine Intelligence}, 4\penalty0 (12):\penalty0 1238--1245, 2022.

\bibitem[Kudela(2023)]{kudela2023evolutionary}
Kudela, J.
\newblock The evolutionary computation methods no one should use.
\newblock \emph{arXiv preprint arXiv:2301.01984}, 2023.

\bibitem[Lange et~al.(2025)Lange, Prasad, Sun, Faldor, Tang, and Ha]{lange2025ai}
Lange, R.~T., Prasad, A., Sun, Q., Faldor, M., Tang, Y., and Ha, D.
\newblock The ai cuda engineer: Agentic cuda kernel discovery, optimization and composition, 2025.

\bibitem[Lewis et~al.(2020)Lewis, Perez, Piktus, Petroni, Karpukhin, Goyal, K{\"u}ttler, Lewis, Yih, Rockt{\"a}schel, et~al.]{lewis2020retrieval}
Lewis, P., Perez, E., Piktus, A., Petroni, F., Karpukhin, V., Goyal, N., K{\"u}ttler, H., Lewis, M., Yih, W.-t., Rockt{\"a}schel, T., et~al.
\newblock Retrieval-augmented generation for knowledge-intensive {NLP} tasks.
\newblock \emph{NeurIPS}, 2020.

\bibitem[Li et~al.(2025{\natexlab{a}})Li, Wang, Bai, Duan, Gao, Hao, and Wen]{li2025formalization1}
Li, C., Wang, Z., Bai, Y., Duan, Y., Gao, Y., Hao, P., and Wen, Z.
\newblock Formalization of algorithms for optimization with block structures.
\newblock \emph{arXiv preprint arXiv:2503.18806}, 2025{\natexlab{a}}.

\bibitem[Li et~al.(2025{\natexlab{b}})Li, Xu, Sun, Zhou, and Wen]{li2025formalization2}
Li, C., Xu, S., Sun, C., Zhou, L., and Wen, Z.
\newblock Formalization of optimality conditions for smooth constrained optimization problems.
\newblock \emph{arXiv preprint arXiv:2503.18821}, 2025{\natexlab{b}}.

\bibitem[Li \& Malik(2017)Li and Malik]{li2017learning}
Li, K. and Malik, J.
\newblock Learning to optimize neural nets.
\newblock \emph{arXiv:1703.00441}, 2017.

\bibitem[Li et~al.(2023)Li, Zhang, and Mak-Hau]{li2023synthesizing}
Li, Q., Zhang, L., and Mak-Hau, V.
\newblock Synthesizing mixed-integer linear programming models from natural language descriptions.
\newblock \emph{arXiv:2311.15271}, 2023.

\bibitem[Li et~al.(2025{\natexlab{c}})Li, Qiao, Wang, Wang, Jin, and Zha]{li2025multi}
Li, W., Qiao, D., Wang, B., Wang, X., Jin, B., and Zha, H.
\newblock Multi-agent credit assignment with pretrained language models.
\newblock In \emph{AISTATS}, 2025{\natexlab{c}}.

\bibitem[Li et~al.(2025{\natexlab{d}})Li, Zhang, Zhang, Zhang, Liu, Yao, Xu, Zheng, Wang, Chen, et~al.]{li2025system}
Li, Z.-Z., Zhang, D., Zhang, M.-L., Zhang, J., Liu, Z., Yao, Y., Xu, H., Zheng, J., Wang, P.-J., Chen, X., et~al.
\newblock From system 1 to system 2: A survey of reasoning large language models.
\newblock \emph{arXiv preprint arXiv:2502.17419}, 2025{\natexlab{d}}.

\bibitem[Liao et~al.(2025)Liao, Wen, Wang, and Zhang]{liao2025marft}
Liao, J., Wen, M., Wang, J., and Zhang, W.
\newblock Marft: Multi-agent reinforcement fine-tuning.
\newblock \emph{arXiv preprint arXiv:2504.16129}, 2025.

\bibitem[Lin et~al.(2024{\natexlab{a}})Lin, Du, Watkins, Hafner, Abbeel, Klein, and Dragan]{lin2024learning}
Lin, J., Du, Y., Watkins, O., Hafner, D., Abbeel, P., Klein, D., and Dragan, A.
\newblock Learning to model the world with language.
\newblock In \emph{ICML}, 2024{\natexlab{a}}.

\bibitem[Lin et~al.(2024{\natexlab{b}})Lin, Tang, Tang, Yang, Chen, Wang, Xiao, Dang, Gan, and Han]{lin2024awq}
Lin, J., Tang, J., Tang, H., Yang, S., Chen, W.-M., Wang, W.-C., Xiao, G., Dang, X., Gan, C., and Han, S.
\newblock {AWQ}: {A}ctivation-aware weight quantization for on-device {LLM} compression and acceleration.
\newblock \emph{GetMobile: Mobile Computing and Communications}, 6:\penalty0 87--100, 2024{\natexlab{b}}.

\bibitem[Lindauer et~al.(2024)Lindauer, Karl, Klier, Moosbauer, Tornede, Mueller, Hutter, Feurer, and Bischl]{lindauer2024position}
Lindauer, M., Karl, F., Klier, A., Moosbauer, J., Tornede, A., Mueller, A., Hutter, F., Feurer, M., and Bischl, B.
\newblock Position: a call to action for a human-centered automl paradigm.
\newblock In \emph{ICML}, 2024.

\bibitem[Liu et~al.(2023)Liu, Tong, Yuan, and Zhang]{liu2023algorithm}
Liu, F., Tong, X., Yuan, M., and Zhang, Q.
\newblock Algorithm evolution using large language model.
\newblock \emph{arXiv:2311.15249}, 2023.

\bibitem[Liu et~al.(2024{\natexlab{a}})Liu, Tong, Yuan, Lin, Luo, Wang, Lu, and Zhang]{EoH}
Liu, F., Tong, X., Yuan, M., Lin, X., Luo, F., Wang, Z., Lu, Z., and Zhang, Q.
\newblock Evolution of heuristics: Towards efficient automatic algorithm design using large language model.
\newblock In \emph{ICML}, 2024{\natexlab{a}}.

\bibitem[Liu et~al.(2024{\natexlab{b}})Liu, Tong, Yuan, Lin, Luo, Wang, Lu, and Zhang]{liu2024example}
Liu, F., Tong, X., Yuan, M., Lin, X., Luo, F., Wang, Z., Lu, Z., and Zhang, Q.
\newblock An example of evolutionary computation+ large language model beating human: {D}esign of efficient guided local search.
\newblock \emph{arXiv:2401.02051}, 2024{\natexlab{b}}.

\bibitem[Liu et~al.(2024{\natexlab{c}})Liu, Yao, Guo, Yang, Lin, Tong, Yuan, Lu, Wang, and Zhang]{liu2024systematic}
Liu, F., Yao, Y., Guo, P., Yang, Z., Lin, X., Tong, X., Yuan, M., Lu, Z., Wang, Z., and Zhang, Q.
\newblock A systematic survey on large language models for algorithm design.
\newblock \emph{arXiv:2410.14716}, 2024{\natexlab{c}}.

\bibitem[Liu et~al.(2018)Liu, Simonyan, and Yang]{liu2018darts}
Liu, H., Simonyan, K., and Yang, Y.
\newblock Darts: Differentiable architecture search.
\newblock \emph{arXiv:1806.09055}, 2018.

\bibitem[Liu et~al.(2024{\natexlab{d}})Liu, Gao, and Li]{liu2024large}
Liu, S., Gao, C., and Li, Y.
\newblock Large language model agent for hyper-parameter optimization.
\newblock \emph{arXiv:2402.01881}, 2024{\natexlab{d}}.

\bibitem[Liu et~al.(2024{\natexlab{e}})Liu, Lou, Jiao, and Zhang]{liu2024position}
Liu, X., Lou, X., Jiao, J., and Zhang, J.
\newblock Position: {F}oundation agents as the paradigm shift for decision making.
\newblock In \emph{ICML}, 2024{\natexlab{e}}.

\bibitem[L{\'o}pez-Ib{\'a}{\~n}ez et~al.(2016)L{\'o}pez-Ib{\'a}{\~n}ez, Dubois-Lacoste, C{\'a}ceres, Birattari, and St{\"u}tzle]{lopez2016irace}
L{\'o}pez-Ib{\'a}{\~n}ez, M., Dubois-Lacoste, J., C{\'a}ceres, L.~P., Birattari, M., and St{\"u}tzle, T.
\newblock The irace package: Iterated racing for automatic algorithm configuration.
\newblock \emph{Operations Research Perspectives}, 3:\penalty0 43--58, 2016.

\bibitem[Luenberger et~al.(1984)Luenberger, Ye, et~al.]{luenberger1984linear}
Luenberger, D.~G., Ye, Y., et~al.
\newblock \emph{Linear and Nonlinear Programming}, volume~2.
\newblock Springer, 1984.

\bibitem[Lv et~al.(2017)Lv, Jiang, and Li]{lv2017learning}
Lv, K., Jiang, S., and Li, J.
\newblock Learning gradient descent: Better generalization and longer horizons.
\newblock In \emph{ICML}, 2017.

\bibitem[Ma et~al.(2024{\natexlab{a}})Ma, Liang, Wang, Huang, Bastani, Jayaraman, Zhu, Fan, and Anandkumar]{ma2024eureka}
Ma, Y.~J., Liang, W., Wang, G., Huang, D.-A., Bastani, O., Jayaraman, D., Zhu, Y., Fan, L., and Anandkumar, A.
\newblock Eureka: {H}uman-level reward design via coding large language models.
\newblock In \emph{ICLR}, 2024{\natexlab{a}}.

\bibitem[Ma et~al.(2024{\natexlab{b}})Ma, Guo, Gong, Zhang, and Tan]{ma2024toward}
Ma, Z., Guo, H., Gong, Y.-J., Zhang, J., and Tan, K.~C.
\newblock Toward automated algorithm design: {A} survey and practical guide to meta-black-box-optimization.
\newblock \emph{arXiv:2411.00625}, 2024{\natexlab{b}}.

\bibitem[Madaan et~al.(2024)Madaan, Tandon, Gupta, Hallinan, Gao, Wiegreffe, Alon, Dziri, Prabhumoye, Yang, et~al.]{madaan2024self}
Madaan, A., Tandon, N., Gupta, P., Hallinan, S., Gao, L., Wiegreffe, S., Alon, U., Dziri, N., Prabhumoye, S., Yang, Y., et~al.
\newblock Self-refine: {I}terative refinement with self-feedback.
\newblock \emph{NeurIPS}, 2024.

\bibitem[Mahammadli(2024)]{mahammadli2024sequential}
Mahammadli, K.
\newblock Sequential large language model-based hyper-parameter optimization.
\newblock \emph{arXiv:2410.20302}, 2024.

\bibitem[Meadows \& Freitas(2022)Meadows and Freitas]{meadows2022survey}
Meadows, J. and Freitas, A.
\newblock A survey in mathematical language processing.
\newblock \emph{arXiv:2205.15231}, 2022.

\bibitem[Metz et~al.(2019)Metz, Maheswaranathan, Nixon, Freeman, and Sohl-Dickstein]{metz2019understanding}
Metz, L., Maheswaranathan, N., Nixon, J., Freeman, D., and Sohl-Dickstein, J.
\newblock Understanding and correcting pathologies in the training of learned optimizers.
\newblock In \emph{ICML}, 2019.

\bibitem[Meyerson et~al.(2024)Meyerson, Nelson, Bradley, Gaier, Moradi, Hoover, and Lehman]{meyerson2024language}
Meyerson, E., Nelson, M.~J., Bradley, H., Gaier, A., Moradi, A., Hoover, A.~K., and Lehman, J.
\newblock Language model crossover: {V}ariation through few-shot prompting.
\newblock \emph{ACM Transactions on Evolutionary Learning}, 4\penalty0 (4):\penalty0 1--40, 2024.

\bibitem[Mirzadeh et~al.(2024)Mirzadeh, Alizadeh, Shahrokhi, Tuzel, Bengio, and Farajtabar]{mirzadeh2024gsm}
Mirzadeh, I., Alizadeh, K., Shahrokhi, H., Tuzel, O., Bengio, S., and Farajtabar, M.
\newblock Gsm-symbolic: Understanding the limitations of mathematical reasoning in large language models.
\newblock \emph{arXiv:2410.05229}, 2024.

\bibitem[Mostajabdaveh et~al.(2024)Mostajabdaveh, Yu, Ramamonjison, Carenini, Zhou, and Zhang]{mostajabdaveh2024optimization}
Mostajabdaveh, M., Yu, T.~T., Ramamonjison, R., Carenini, G., Zhou, Z., and Zhang, Y.
\newblock Optimization modeling and verification from problem specifications using a multi-agent multi-stage {LLM} framework.
\newblock \emph{INFOR: Information Systems and Operational Research}, 62\penalty0 (4):\penalty0 599--617, 2024.

\bibitem[Moura \& Ullrich(2021)Moura and Ullrich]{moura2021lean}
Moura, L.~d. and Ullrich, S.
\newblock The {L}ean 4 theorem prover and programming language.
\newblock In \emph{CADE}, 2021.

\bibitem[Nelder \& Mead(1965)Nelder and Mead]{nelder1965simplex}
Nelder, J.~A. and Mead, R.
\newblock A simplex method for function minimization.
\newblock \emph{The Computer Journal}, 7\penalty0 (4):\penalty0 308--313, 1965.

\bibitem[Ning et~al.(2023)Ning, Liu, Qin, Xiao, Xue, Huang, Liu, Chen, and Wu]{ning2023novel}
Ning, Y., Liu, J., Qin, L., Xiao, T., Xue, S., Huang, Z., Liu, Q., Chen, E., and Wu, J.
\newblock A novel approach for auto-formulation of optimization problems.
\newblock \emph{arXiv:2302.04643}, 2023.

\bibitem[Oliveira et~al.(2016)Oliveira, Lima, and Montevechi]{oliveira2016perspectives}
Oliveira, J.~B., Lima, R.~S., and Montevechi, J. A.~B.
\newblock Perspectives and relationships in supply chain simulation: {A} systematic literature review.
\newblock \emph{Simulation Modelling Practice and Theory}, 62:\penalty0 166--191, 2016.

\bibitem[{OpenAI}(2024)]{openai2024o1}
{OpenAI}.
\newblock Learning to reason with llms.
\newblock \url{https://openai.com/index/learning-to-reason-with-llms/}, 2024.

\bibitem[OpenAI(2024)]{rft}
OpenAI.
\newblock Reinforcement fine-tuning, 2024.
\newblock URL \url{https://www.youtube.com/watch?v=yCIYS9fx56U}.

\bibitem[Ormerod \& Kiossis(1997)Ormerod and Kiossis]{ormerod1997or}
Ormerod, R. and Kiossis, I.
\newblock {OR}/{MS} publications: {E}xtension of the analysis of us flagship journals to the united kingdom.
\newblock \emph{Operations Research}, 45\penalty0 (2):\penalty0 178--187, 1997.

\bibitem[Ormerod(2002)]{ormerod2002nature}
Ormerod, R.~J.
\newblock On the nature of {OR}: {T}aking stock.
\newblock \emph{Journal of the Operational Research Society}, 53\penalty0 (5):\penalty0 475--491, 2002.

\bibitem[Park et~al.(2023)Park, O'Brien, Cai, Morris, Liang, and Bernstein]{park2023generative}
Park, J.~S., O'Brien, J., Cai, C.~J., Morris, M.~R., Liang, P., and Bernstein, M.~S.
\newblock Generative agents: {I}nteractive simulacra of human behavior.
\newblock In \emph{UIST}, 2023.

\bibitem[Pham et~al.(2018)Pham, Guan, Zoph, Le, and Dean]{pham2018efficient}
Pham, H., Guan, M., Zoph, B., Le, Q., and Dean, J.
\newblock Efficient neural architecture search via parameters sharing.
\newblock In \emph{ICML}, pp.\  4095--4104, 2018.

\bibitem[Pietri \& Sakellariou(2016)Pietri and Sakellariou]{pietri2016mapping}
Pietri, I. and Sakellariou, R.
\newblock Mapping virtual machines onto physical machines in cloud computing: {A} survey.
\newblock \emph{ACM Computing Surveys (CSUR)}, 49\penalty0 (3):\penalty0 1--30, 2016.

\bibitem[Pluhacek et~al.(2023)Pluhacek, Kazikova, Kadavy, Viktorin, and Senkerik]{pluhacek2023leveraging}
Pluhacek, M., Kazikova, A., Kadavy, T., Viktorin, A., and Senkerik, R.
\newblock Leveraging large language models for the generation of novel metaheuristic optimization algorithms.
\newblock In \emph{GECCO}, 2023.

\bibitem[Ramamonjison et~al.(2022)Ramamonjison, Li, Yu, He, Rengan, Banitalebi-Dehkordi, Zhou, and Zhang]{ramamonjison2022augmenting}
Ramamonjison, R., Li, H., Yu, T.~T., He, S., Rengan, V., Banitalebi-Dehkordi, A., Zhou, Z., and Zhang, Y.
\newblock Augmenting operations research with auto-formulation of optimization models from problem descriptions.
\newblock \emph{arXiv:2209.15565}, 2022.

\bibitem[Ramamonjison et~al.(2023)Ramamonjison, Yu, Li, Li, Carenini, Ghaddar, He, Mostajabdaveh, Banitalebi-Dehkordi, Zhou, et~al.]{ramamonjison2023nl4opt}
Ramamonjison, R., Yu, T., Li, R., Li, H., Carenini, G., Ghaddar, B., He, S., Mostajabdaveh, M., Banitalebi-Dehkordi, A., Zhou, Z., et~al.
\newblock {NL4OPT} competition: {F}ormulating optimization problems based on their natural language descriptions.
\newblock In \emph{NeurIPS 2022 Competition Track}, 2023.

\bibitem[Real et~al.(2019)Real, Aggarwal, Huang, and Le]{real2019regularized}
Real, E., Aggarwal, A., Huang, Y., and Le, Q.~V.
\newblock Regularized evolution for image classifier architecture search.
\newblock In \emph{AAAI}, 2019.

\bibitem[Reisman \& Kirschnick(1994)Reisman and Kirschnick]{reisman1994devolution}
Reisman, A. and Kirschnick, F.
\newblock The devolution of {OR}/{MS}: Implications from a statistical content analysis of papers in flagship journals.
\newblock \emph{Operations Research}, 42\penalty0 (4):\penalty0 577--588, 1994.

\bibitem[Romera{-}Paredes et~al.(2024)Romera{-}Paredes, Barekatain, Novikov, Balog, Kumar, Dupont, Ruiz, Ellenberg, Wang, Fawzi, Kohli, and Fawzi]{Funsearch}
Romera{-}Paredes, B., Barekatain, M., Novikov, A., Balog, M., Kumar, M.~P., Dupont, E., Ruiz, F., Ellenberg, J.~S., Wang, P., Fawzi, O., Kohli, P., and Fawzi, A.
\newblock Mathematical discoveries from program search with large language models.
\newblock \emph{Nature}, 625\penalty0 (7995):\penalty0 468--475, 2024.

\bibitem[Romero~Morales(2000)]{romero2000optimization}
Romero~Morales, M.~D.
\newblock \emph{Optimization Problems in Supply Chain Management}.
\newblock PhD thesis, TRAIL Research School, The Netherlands, 2000.
\newblock ERIM PhD series Research in Management nr. 3.

\bibitem[Rossit et~al.(2019)Rossit, Vigo, Tohm{\'e}, and Frutos]{rossit2019visual}
Rossit, D.~G., Vigo, D., Tohm{\'e}, F., and Frutos, M.
\newblock Visual attractiveness in routing problems: {A} review.
\newblock \emph{Computers \& Operations Research}, 103:\penalty0 13--34, 2019.

\bibitem[Roziere et~al.(2023)Roziere, Gehring, Gloeckle, Sootla, Gat, Tan, Adi, Liu, Sauvestre, Remez, et~al.]{roziere2023code}
Roziere, B., Gehring, J., Gloeckle, F., Sootla, S., Gat, I., Tan, X.~E., Adi, Y., Liu, J., Sauvestre, R., Remez, T., et~al.
\newblock Code {L}lama: {O}pen foundation models for code.
\newblock \emph{arXiv:2308.12950}, 2023.

\bibitem[Ruder(2016)]{ruder2016overview}
Ruder, S.
\newblock An overview of gradient descent optimization algorithms.
\newblock \emph{arXiv:1609.04747}, 2016.

\bibitem[Schick \& Sch{\"u}tze(2020)Schick and Sch{\"u}tze]{schick2020s}
Schick, T. and Sch{\"u}tze, H.
\newblock It's not just size that matters: Small language models are also few-shot learners.
\newblock \emph{arXiv:2009.07118}, 2020.

\bibitem[Schick et~al.(2023)Schick, Dwivedi-Yu, Dess{\`\i}, Raileanu, Lomeli, Hambro, Zettlemoyer, Cancedda, and Scialom]{schick2023toolformer}
Schick, T., Dwivedi-Yu, J., Dess{\`\i}, R., Raileanu, R., Lomeli, M., Hambro, E., Zettlemoyer, L., Cancedda, N., and Scialom, T.
\newblock Toolformer: {L}anguage models can teach themselves to use tools.
\newblock \emph{NeurIPS}, 2023.

\bibitem[Schmidt et~al.(2007)Schmidt, Fung, and Rosales]{schmidt2007fast}
Schmidt, M., Fung, G., and Rosales, R.
\newblock Fast optimization methods for $\ell_1$ regularization: A comparative study and two new approaches.
\newblock In \emph{ECML}, 2007.

\bibitem[Shen et~al.(2024)Shen, Song, Tan, Li, Lu, and Zhuang]{shen2024hugginggpt}
Shen, Y., Song, K., Tan, X., Li, D., Lu, W., and Zhuang, Y.
\newblock Hugginggpt: {S}olving ai tasks with chatgpt and its friends in hugging face.
\newblock \emph{NeurIPS}, 2024.

\bibitem[Sheng et~al.(2022)Sheng, Hu, Zhou, Zhu, Jin, Wang, and Wang]{sheng2022learning}
Sheng, J., Hu, Y., Zhou, W., Zhu, L., Jin, B., Wang, J., and Wang, X.
\newblock Learning to schedule multi-{NUMA} virtual machines via reinforcement learning.
\newblock \emph{Pattern Recognition}, 121:\penalty0 108254, 2022.

\bibitem[Shin(2021)]{shin2021effects}
Shin, D.
\newblock The effects of explainability and causability on perception, trust, and acceptance: {I}mplications for explainable {AI}.
\newblock \emph{International journal of human-computer studies}, 146:\penalty0 102551, 2021.

\bibitem[Shinn et~al.(2024)Shinn, Cassano, Gopinath, Narasimhan, and Yao]{shinn2024reflexion}
Shinn, N., Cassano, F., Gopinath, A., Narasimhan, K., and Yao, S.
\newblock Reflexion: Language agents with verbal reinforcement learning.
\newblock \emph{NeurIPS}, 2024.

\bibitem[Silver \& Sutton(2025)Silver and Sutton]{silver2025experience}
Silver, D. and Sutton, R.~S.
\newblock Welcome to the era of experience.
\newblock \emph{Preprint of a chapter in Designing an Intelligence}, 2025.

\bibitem[Simchi-Levi(2014)]{simchi2014om}
Simchi-Levi, D.
\newblock {OM} forum—{OM} research: From problem-driven to data-driven research.
\newblock \emph{Manufacturing \& Service Operations Management}, 16\penalty0 (1):\penalty0 2--10, 2014.

\bibitem[Snoek et~al.(2012)Snoek, Larochelle, and Adams]{snoek2012practical}
Snoek, J., Larochelle, H., and Adams, R.~P.
\newblock Practical bayesian optimization of machine learning algorithms.
\newblock \emph{NeurIPS}, 25, 2012.

\bibitem[Song et~al.(2024)Song, Yang, and Anandkumar]{song2024towards}
Song, P., Yang, K., and Anandkumar, A.
\newblock Towards large language models as copilots for theorem proving in lean.
\newblock \emph{arXiv:2404.12534}, 2024.

\bibitem[Sun et~al.(2022)Sun, Zhang, Hu, and Van~Mieghem]{sun2022predicting}
Sun, J., Zhang, D.~J., Hu, H., and Van~Mieghem, J.~A.
\newblock Predicting human discretion to adjust algorithmic prescription: {A} large-scale field experiment in warehouse operations.
\newblock \emph{Management Science}, 68\penalty0 (2):\penalty0 846--865, 2022.

\bibitem[Sun et~al.(2019)Sun, Cao, Zhu, and Zhao]{sun2019survey}
Sun, S., Cao, Z., Zhu, H., and Zhao, J.
\newblock A survey of optimization methods from a machine learning perspective.
\newblock \emph{IEEE Transactions on Cybernetics}, 50\penalty0 (8):\penalty0 3668--3681, 2019.

\bibitem[Tang et~al.(2024)Tang, Huang, Zheng, Hu, Wang, Ge, and Wang]{tang2024orlm}
Tang, Z., Huang, C., Zheng, X., Hu, S., Wang, Z., Ge, D., and Wang, B.
\newblock {ORLM}: Training large language models for optimization modeling.
\newblock \emph{arXiv:2405.17743}, 2024.

\bibitem[Trirat et~al.(2024)Trirat, Jeong, and Hwang]{trirat2024automl}
Trirat, P., Jeong, W., and Hwang, S.~J.
\newblock Automl-agent: {A} multi-agent llm framework for full-pipeline automl.
\newblock \emph{arXiv:2410.02958}, 2024.

\bibitem[Tsouros et~al.(2023)Tsouros, Verhaeghe, Kad{\i}o{\u{g}}lu, and Guns]{tsouros2023holy}
Tsouros, D., Verhaeghe, H., Kad{\i}o{\u{g}}lu, S., and Guns, T.
\newblock Holy {G}rail 2.0: {F}rom natural language to constraint models.
\newblock \emph{arXiv:2308.01589}, 2023.

\bibitem[Valmeekam et~al.(2022)Valmeekam, Olmo, Sreedharan, and Kambhampati]{valmeekam2022large}
Valmeekam, K., Olmo, A., Sreedharan, S., and Kambhampati, S.
\newblock Large language models still can't plan (a benchmark for {LLMs} on planning and reasoning about change).
\newblock In \emph{NeurIPS 2022 Foundation Models for Decision Making Workshop}, 2022.

\bibitem[Valmeekam et~al.(2024)Valmeekam, Stechly, and Kambhampati]{valmeekam2024llms}
Valmeekam, K., Stechly, K., and Kambhampati, S.
\newblock {LLMs} still can't plan; {C}an lrms? {A} preliminary evaluation of {OpenAI's} o1 on planbench.
\newblock \emph{arXiv:2409.13373}, 2024.

\bibitem[Van~Donselaar et~al.(2010)Van~Donselaar, Gaur, Van~Woensel, Broekmeulen, and Fransoo]{van2010ordering}
Van~Donselaar, K.~H., Gaur, V., Van~Woensel, T., Broekmeulen, R.~A., and Fransoo, J.~C.
\newblock Ordering behavior in retail stores and implications for automated replenishment.
\newblock \emph{Management Science}, 56\penalty0 (5):\penalty0 766--784, 2010.

\bibitem[Van~Nguyen et~al.(2024)Van~Nguyen, Shen, Aponte, Xia, Basu, Hu, Chen, Parmar, Kunapuli, Barrow, et~al.]{van2024survey}
Van~Nguyen, C., Shen, X., Aponte, R., Xia, Y., Basu, S., Hu, Z., Chen, J., Parmar, M., Kunapuli, S., Barrow, J., et~al.
\newblock A survey of small language models.
\newblock \emph{arXiv:2410.20011}, 2024.

\bibitem[Velasco et~al.(2024)Velasco, Guerrero, and Hospitaler]{velasco2024literature}
Velasco, L., Guerrero, H., and Hospitaler, A.
\newblock A literature review and critical analysis of metaheuristics recently developed.
\newblock \emph{Archives of Computational Methods in Engineering}, 31\penalty0 (1):\penalty0 125--146, 2024.

\bibitem[Wang et~al.(2024)Wang, Xie, Jiang, Mandlekar, Xiao, Zhu, Fan, and Anandkumar]{wang2024voyager}
Wang, G., Xie, Y., Jiang, Y., Mandlekar, A., Xiao, C., Zhu, Y., Fan, L., and Anandkumar, A.
\newblock Voyager: {A}n open-ended embodied agent with large language models.
\newblock \emph{Transactions on Machine Learning Research}, 2024.
\newblock ISSN 2835-8856.

\bibitem[Wang et~al.(2023)Wang, Cai, Chen, Liu, Ma, and Liang]{wang2023describe}
Wang, Z., Cai, S., Chen, G., Liu, A., Ma, X.~S., and Liang, Y.
\newblock Describe, explain, plan and select: {I}nteractive planning with llms enables open-world multi-task agents.
\newblock \emph{NeurIPS}, 2023.

\bibitem[Wasserkrug et~al.(2024)Wasserkrug, Boussioux, Hertog, Mirzazadeh, Birbil, Kurtz, and Maragno]{wasserkrug2024large}
Wasserkrug, S., Boussioux, L., Hertog, D.~d., Mirzazadeh, F., Birbil, I., Kurtz, J., and Maragno, D.
\newblock From large language models and optimization to decision optimization copilot: A research manifesto.
\newblock \emph{arXiv:2402.16269}, 2024.

\bibitem[Wei et~al.(2022)Wei, Wang, Schuurmans, Bosma, Xia, Chi, Le, Zhou, et~al.]{wei2022chain}
Wei, J., Wang, X., Schuurmans, D., Bosma, M., Xia, F., Chi, E., Le, Q.~V., Zhou, D., et~al.
\newblock Chain-of-thought prompting elicits reasoning in large language models.
\newblock \emph{NeurIPS}, 2022.

\bibitem[Wohlberg(2017)]{wohlberg2017admm}
Wohlberg, B.
\newblock {ADMM} penalty parameter selection by residual balancing.
\newblock \emph{arXiv:1704.06209}, 2017.

\bibitem[Woldesenbet \& Yen(2009)Woldesenbet and Yen]{woldesenbet2009dynamic}
Woldesenbet, Y.~G. and Yen, G.~G.
\newblock Dynamic evolutionary algorithm with variable relocation.
\newblock \emph{IEEE Transactions on Evolutionary Computation}, 13\penalty0 (3):\penalty0 500--513, 2009.

\bibitem[Wolke et~al.(2015)Wolke, Tsend-Ayush, Pfeiffer, and Bichler]{wolke2015more}
Wolke, A., Tsend-Ayush, B., Pfeiffer, C., and Bichler, M.
\newblock More than bin packing: {D}ynamic resource allocation strategies in cloud data centers.
\newblock \emph{Information Systems}, 52:\penalty0 83--95, 2015.

\bibitem[Wolpert \& Macready(1997)Wolpert and Macready]{wolpert1997no}
Wolpert, D.~H. and Macready, W.~G.
\newblock No free lunch theorems for optimization.
\newblock \emph{IEEE Transactions on Evolutionary Computation}, 1\penalty0 (1):\penalty0 67--82, 1997.

\bibitem[Wu et~al.(2024)Wu, Wu, Wu, Feng, and Tan]{wu2024evolutionary}
Wu, X., Wu, S.-h., Wu, J., Feng, L., and Tan, K.~C.
\newblock Evolutionary computation in the era of large language model: {S}urvey and roadmap.
\newblock \emph{arXiv:2401.10034}, 2024.

\bibitem[Wu et~al.(2022)Wu, Jiang, Li, Rabe, Staats, Jamnik, and Szegedy]{wu2022autoformalization}
Wu, Y., Jiang, A.~Q., Li, W., Rabe, M., Staats, C., Jamnik, M., and Szegedy, C.
\newblock Autoformalization with large language models.
\newblock \emph{NeurIPS}, 2022.

\bibitem[Xi et~al.(2023)Xi, Chen, Guo, He, Ding, Hong, Zhang, Wang, Jin, Zhou, et~al.]{xi2023rise}
Xi, Z., Chen, W., Guo, X., He, W., Ding, Y., Hong, B., Zhang, M., Wang, J., Jin, S., Zhou, E., et~al.
\newblock The rise and potential of large language model based agents: {A} survey.
\newblock \emph{arXiv:2309.07864}, 2023.

\bibitem[Xu et~al.(2025)Xu, Hao, Zong, Wang, Zhang, Wang, Lan, Gong, Ouyang, Meng, et~al.]{xu2025towards}
Xu, F., Hao, Q., Zong, Z., Wang, J., Zhang, Y., Wang, J., Lan, X., Gong, J., Ouyang, T., Meng, F., et~al.
\newblock Towards large reasoning models: A survey of reinforced reasoning with large language models.
\newblock \emph{arXiv preprint arXiv:2501.09686}, 2025.

\bibitem[Xu et~al.(2024)Xu, Li, Tao, Shen, Cheng, Li, Xu, Tao, and Zhou]{xu2024survey}
Xu, X., Li, M., Tao, C., Shen, T., Cheng, R., Li, J., Xu, C., Tao, D., and Zhou, T.
\newblock A survey on knowledge distillation of large language models.
\newblock \emph{arXiv:2402.13116}, 2024.

\bibitem[Xu et~al.(2017{\natexlab{a}})Xu, Liu, Lin, and Yang]{xu2017admm}
Xu, Y., Liu, M., Lin, Q., and Yang, T.
\newblock Admm without a fixed penalty parameter: Faster convergence with new adaptive penalization.
\newblock \emph{NeurIPS}, 2017{\natexlab{a}}.

\bibitem[Xu et~al.(2017{\natexlab{b}})Xu, Figueiredo, and Goldstein]{xu2017adaptive}
Xu, Z., Figueiredo, M., and Goldstein, T.
\newblock Adaptive admm with spectral penalty parameter selection.
\newblock In \emph{AISTATS}, 2017{\natexlab{b}}.

\bibitem[Yadkori et~al.(2024)Yadkori, Kuzborskij, Gy{\"o}rgy, and Szepesv{\'a}ri]{yadkori2024believe}
Yadkori, Y.~A., Kuzborskij, I., Gy{\"o}rgy, A., and Szepesv{\'a}ri, C.
\newblock To believe or not to believe your llm.
\newblock \emph{arXiv:2406.02543}, 2024.

\bibitem[Yampolskiy(2018)]{yampolskiy2018we}
Yampolskiy, R.~V.
\newblock Why we do not evolve software? {A}nalysis of evolutionary algorithms.
\newblock \emph{Evolutionary Bioinformatics}, 14:\penalty0 1176934318815906, 2018.

\bibitem[Yang et~al.(2024)Yang, Wang, Lu, Liu, Le, Zhou, and Chen]{yang2024large}
Yang, C., Wang, X., Lu, Y., Liu, H., Le, Q.~V., Zhou, D., and Chen, X.
\newblock Large language models as optimizers.
\newblock In \emph{ICLR}, 2024.

\bibitem[Yang et~al.(2022)Yang, Guan, Jia, Yu, Xu, and Spanos]{yang2022survey}
Yang, Y., Guan, X., Jia, Q.-S., Yu, L., Xu, B., and Spanos, C.~J.
\newblock A survey of {ADMM} variants for distributed optimization: {P}roblems, algorithms and features.
\newblock \emph{arXiv:2208.03700}, 2022.

\bibitem[Yao et~al.(2018)Yao, Wang, Chen, Dai, Li, Tu, Yang, and Yu]{yao2018taking}
Yao, Q., Wang, M., Chen, Y., Dai, W., Li, Y.-F., Tu, W.-W., Yang, Q., and Yu, Y.
\newblock Taking human out of learning applications: {A} survey on automated machine learning.
\newblock \emph{arXiv:1810.13306}, 2018.

\bibitem[Yao et~al.(2022)Yao, Zhao, Yu, Du, Shafran, Narasimhan, and Cao]{yao2022react}
Yao, S., Zhao, J., Yu, D., Du, N., Shafran, I., Narasimhan, K., and Cao, Y.
\newblock React: {S}ynergizing reasoning and acting in language models.
\newblock \emph{arXiv:2210.03629}, 2022.

\bibitem[Yao et~al.(2024{\natexlab{a}})Yao, Liu, Lin, Lu, Wang, and Zhang]{yao2024multi}
Yao, S., Liu, F., Lin, X., Lu, Z., Wang, Z., and Zhang, Q.
\newblock Multi-objective evolution of heuristic using large language model.
\newblock \emph{arXiv:2409.16867}, 2024{\natexlab{a}}.

\bibitem[Yao et~al.(2024{\natexlab{b}})Yao, Yu, Zhao, Shafran, Griffiths, Cao, and Narasimhan]{yao2024tree}
Yao, S., Yu, D., Zhao, J., Shafran, I., Griffiths, T., Cao, Y., and Narasimhan, K.
\newblock Tree of thoughts: {D}eliberate problem solving with large language models.
\newblock \emph{NeurIPS}, 2024{\natexlab{b}}.

\bibitem[Ye et~al.(2024)Ye, Wang, Cao, Berto, Hua, Kim, Park, and Song]{ye2024reevo}
Ye, H., Wang, J., Cao, Z., Berto, F., Hua, C., Kim, H., Park, J., and Song, G.
\newblock Re{E}vo: {L}arge language models as hyper-heuristics with reflective evolution.
\newblock \emph{arXiv:2402.01145}, 2024.

\bibitem[Yu \& Gen(2010)Yu and Gen]{yu2010introduction}
Yu, X. and Gen, M.
\newblock \emph{Introduction to evolutionary algorithms}.
\newblock Springer Science \& Business Media, 2010.

\bibitem[Zhang et~al.(2019)Zhang, Wang, Zhang, Dai, and Shen]{zhang2019gap}
Zhang, D., Wang, L., Zhang, L., Dai, B.~T., and Shen, H.~T.
\newblock The gap of semantic parsing: A survey on automatic math word problem solvers.
\newblock \emph{IEEE Transactions on Pattern Analysis and Machine Intelligence}, 42\penalty0 (9):\penalty0 2287--2305, 2019.

\bibitem[Zhang et~al.(2023{\natexlab{a}})Zhang, Desai, Bae, Lorraine, and Ba]{zhang2023using}
Zhang, M., Desai, N., Bae, J., Lorraine, J., and Ba, J.
\newblock Using large language models for hyperparameter optimization.
\newblock In \emph{NeurIPS 2023 Foundation Models for Decision Making Workshop}, 2023{\natexlab{a}}.

\bibitem[Zhang et~al.(2024)Zhang, Liu, Lin, Wang, Lu, and Zhang]{zhang2024understanding}
Zhang, R., Liu, F., Lin, X., Wang, Z., Lu, Z., and Zhang, Q.
\newblock Understanding the importance of evolutionary search in automated heuristic design with large language models.
\newblock In \emph{PPSN}, 2024.

\bibitem[Zhang et~al.(2023{\natexlab{b}})Zhang, Gong, Wu, Liu, and Zhou]{zhang2023automl}
Zhang, S., Gong, C., Wu, L., Liu, X., and Zhou, M.
\newblock {AutoML-GPT}: {A}utomatic machine learning with {GPT}.
\newblock \emph{arXiv:2305.02499}, 2023{\natexlab{b}}.

\bibitem[Zheng et~al.(2025{\natexlab{a}})Zheng, Qiu, Shi, and Ma]{zheng2025towards}
Zheng, J., Qiu, S., Shi, C., and Ma, Q.
\newblock Towards lifelong learning of large language models: A survey.
\newblock \emph{ACM Computing Surveys}, 57\penalty0 (8):\penalty0 1--35, 2025{\natexlab{a}}.

\bibitem[Zheng et~al.(2025{\natexlab{b}})Zheng, Shi, Cai, Li, Zhang, Li, Yu, and Ma]{zheng2025lifelong}
Zheng, J., Shi, C., Cai, X., Li, Q., Zhang, D., Li, C., Yu, D., and Ma, Q.
\newblock Lifelong learning of large language model based agents: A roadmap.
\newblock \emph{arXiv preprint arXiv:2501.07278}, 2025{\natexlab{b}}.

\bibitem[Zhong et~al.(2024)Zhong, Xu, Zhang, and Yu]{zhong2024leveraging}
Zhong, R., Xu, Y., Zhang, C., and Yu, J.
\newblock Leveraging large language model to generate a novel metaheuristic algorithm with {CRISPE} framework.
\newblock \emph{Cluster Computing}, 27\penalty0 (10):\penalty0 13835--13869, 2024.

\bibitem[Zhou et~al.(2011)Zhou, Qu, Li, Zhao, Suganthan, and Zhang]{zhou2011multiobjective}
Zhou, A., Qu, B.-Y., Li, H., Zhao, S.-Z., Suganthan, P.~N., and Zhang, Q.
\newblock Multiobjective evolutionary algorithms: A survey of the state of the art.
\newblock \emph{Swarm and evolutionary computation}, 1\penalty0 (1):\penalty0 32--49, 2011.

\bibitem[Zhou et~al.(2025)Zhou, Tie, Zhang, Wang, Zuo, Wu, Chu, Zhou, Sun, and Gong]{zhou2025large}
Zhou, X., Tie, G., Zhang, G., Wang, W., Zuo, Z., Wu, D., Chu, D., Zhou, P., Sun, L., and Gong, N.~Z.
\newblock Large reasoning models in agent scenarios: Exploring the necessity of reasoning capabilities.
\newblock \emph{arXiv preprint arXiv:2503.11074}, 2025.

\end{thebibliography}
\bibliographystyle{icml2025}


\end{document}